\documentclass[12pt, reqno]{amsart}

\usepackage{amsmath, amssymb, amscd, amsthm, comment, amsxtra}
\usepackage{graphicx}
\usepackage{subfigure}

\usepackage[pdftex, linktocpage=true]{hyperref}
\usepackage{euscript}

\usepackage[format=plain,labelfont=bf,up,width=.99\textwidth]{caption}

\newcommand{\field}[1]{\mathbb{#1}} 
\newtheorem{thm}{Theorem}[section] 
\newtheorem*{mainprop*}{Main Proposition}
\newtheorem*{thmA*}{Theorem A}
\newtheorem*{thmB*}{Theorem B}
\newtheorem*{KPr*}{Theorem (Koebe Priniciple)}
\newtheorem*{OKPr*}{Theorem (One-sided Koebe Priniciple)}

\newtheorem{lemma}[thm]{Lemma} 
\newtheorem{prop}[thm]{Proposition}
\newtheorem{cor}[thm]{Corollary}
\newtheorem{definition}[thm]{Definition}
 
\newtheorem{remark}[thm]{Remark}

\def\eps{\epsilon}
\def\vareps{\varepsilon}

\def\cA{\mathcal{A}}
\def\cB{\mathcal{B}}

\def\cD{\mathcal{D}}

\def\cK{\mathcal{K}}

\def\cL{\mathcal{L}}

\def\cM{\mathcal{M}}

\def\cR{\mathcal{R}}

\def\cP{\mathcal{P}}
\def\cS{\mathcal{S}}

\def\cY{\mathcal{Y}}

\date{\today}

\begin{document}

\title{Renormalization for Lorenz maps of  monotone combinatorial types}
\author{Denis Gaidashev}
\address{Department of Mathematics, Uppsala University, Uppsala, Sweden.}

\begin{abstract}
Lorenz maps are maps of the unit interval with one critical point of order $\rho>1$, and a discontinuity at that point. They appear as return maps of sections of the geometric Lorenz flow.

We construct real {\it a priori} bounds for renormalizable Lorenz maps with certain monotone combinatorics and a sufficiently flat critical point, and use these bounds to show existence of periodic points of renormalization, as well as existence of Cantor attractors for dynamics of infinitely renormalizable Lorenz maps.
\end{abstract}
\maketitle

\tableofcontents

\setcounter{page}{1}


\section{Introduction}

E. N. Lorenz in \cite{Lor} demonstrated  numerically the existence of certain three-dimensional flows that have a complicated behaviour. The {\it Lorenz flow} has a saddle fixed point with a one-dimensional unstable manifold and an infinite set of periodic orbits whose closure constitutes a global attractor of the flow.

As it is often done in dynamics, one can attempt to understand the behaviour of a three-dimensional flow by looking at the first return map to an appropriately chosen two-dimensional section. In the case of the Lorenz flow, it is convenient to choose the section as a plane transversal to the local stable manifold, and, therefore,  intersecting it along a curve $\gamma$. The first return map is discontinuous at $\gamma$.

The {\it geometric Lorenz flow} has been introduced in \cite{Wil}: a Lorenz flow with an extra condition that the return map preserves a one-dimensional foliation in the section, and contracts distances between points in the leafs of this foliation at a geometric rate. Since the return map is contracting in the leafs, its dynamics is asymptotically one-dimensional, and can be understood in terms of a map acting on the space of leafs (an interval). This interval map has a discontinuity  at the point of the interval corresponding to $\gamma$, and is commonly called  the {\it Lorenz map}. 

\vspace{3mm}

\hfill\begin{minipage}{\dimexpr\textwidth-6mm}
{\it In this paper, we construct renormalization {real a priori bounds} for a class of Lorenz maps of bounded monotone combinatorics (Main Proposition and Theorem B, Section $\ref{statement_results}$), and demonstrate existence of {uniquely ergodic Cantor attractors} for infinitely renormalizable Lorenz maps in this class (Theorem A, Section $\ref{statement_results}$).}
\end{minipage}

\vspace{3mm}

Similar results have been obtained for another class of Lorenz maps of monotone combinatorics in \cite{Win2}. In the present paper we extend the approach of \cite{Win2} to a rather different combinatorial class of Lorenz maps (see explanation after theorem B for details).


\subsection{Background and definitions}

We will start by defining what is known as the standard Lorenz family. Our work is a continuation of the study started in \cite{Win2}, and we will, therefore, make a conscientious effort to use the notation of \cite{Win2} so that it would be easier for the reader to compare the approach of this paper  with that of \cite{Win2}.
\begin{definition} \label{sLf}
Let $u \in [0,1]$, $v \in [0,1]$, $c \in (0,1)$ and $\rho > 0$. The standard Lorenz family $(u,v,c) \mapsto Q(x)$ is the family of maps $Q:[0,1] \setminus \{ c \} \mapsto [0,1]$ defined as 
$$Q(x)=\left\{
\begin{tabular}{l l}
$u \left(1-\left(c-x \over c   \right)^\rho \right)$, & $x  \in  [0,c)$, \\  
$1+ v \left(-1+\left(x-c \over 1-c   \right)^\rho \right)$,& $x  \in  (c,1]$,   
\end{tabular}
\right.$$
\end{definition}

\begin{remark}
In the definition above, $u$ is the length of $Q([0,c))$, $v$ is that of $Q((c,1])$, while $u$ and $1-v$ are the critical values. To emphasize that a critical point $c$ corresponds to a map $f$, we will use the notation $c(f)$. The difference $1-c$ will be denoted as $\mu$:
$$\mu \equiv 1-c.$$
\end{remark}

More generally,
\begin{definition} \label{Lorenz_family}
A $C^k$-Lorenz map $f:[0,1] \setminus \{c\} \mapsto [0,1]$ is defined as a triple $(\phi,\psi,Q)$, 
$$f(x)=\left\{
\begin{tabular}{l l}
$\phi(Q(x))$,& $x  \in  [0,c)$, \\  
$\psi(Q(x))$,& $x  \in  (c,1]$,   
\end{tabular}
\right.$$
where $\phi$ and $\psi$ are $C^k$ orientation preserving diffeomorphisms of $[0,1]$ (this space will be denoted by $\cD^k$). 
\end{definition}

\begin{remark}
Notice, two different $C^k$-Lorenz maps, that is two different triples $(\phi,\psi,Q)$ and $(\tilde \phi,\tilde \psi,\tilde Q)$ define one and the same map on $[0,1]$ if $c=\tilde c$ and  $\tilde \phi \arrowvert_{[0,c)}=\phi \circ u/\tilde u \arrowvert_{[0,c)}$ and $\tilde \psi \arrowvert_{[0,c)}=\psi \circ (v (x-1)/\tilde v +1) \arrowvert_{[0,c)}$. We will, however, think of two triples $(\phi,\psi,Q)$ and $(\tilde \phi,\tilde \psi, \tilde Q)$ as distinct. The ambiguity disappears in the class $C^\omega$, since, e.g.,  $\tilde \phi \arrowvert_{[0,c)}=\phi \circ u/\tilde u \arrowvert_{[0,c)} \implies  \tilde \phi \arrowvert_{[0,1]}=\phi \circ u/\tilde u \arrowvert_{[0,1]} \implies \tilde \phi(1) \ne 1$.
\end{remark}

We will refer to the diffeomorphisms $\phi$ and $\psi$ as {\it coefficients} or {\it factors}  of the Lorenz map.

The set of $C^k$-Lorenz maps will be denoted $\cL^k$. Since a Lorenz map $(\ref{Lorenz_family})$ can be identified with a quintuple $(u,v,c,\phi,\psi)$, the space $\cL^k$ is isomorphic to $[0,1]^2 \times (0,1) \times \cD^k \times \cD^k$.  $\cL^S \subset \cL^3$ will denote the subset of maps with  negative Schwarzian derivative $S_f$,
\begin{equation}\label{Schw}
S_f(x)={f'''(x) \over f'(x)  }-{3  \over 2} \left({f''(x) \over f'(x) } \right)^2
\end{equation}
 The notation $| \cdot |_k$ will be used for the $C^k$-norm.  The subsets  of $\cD^3$ of diffeomorphisms with a negative Schwarzian will be denoted $\cD^S$.

Guckenheimer and Williams have proved in \cite{GW} that there is an open set of three-dimensional vector fields, that generate a geometric Lorenz flow with a smooth Lorenz map of $\rho<1$. However, one can use the arguments of \cite{GW} to construct open sets of vector fields with Lorenz maps of $\rho \ge 1$.  Similarly to the unimodal family, Lorenz maps with $\rho > 1$ have  richer dynamics that combines contraction with  expansion.

\begin{definition} \label{branch}
A branch of $f^n$ is a maximal closed interval $J$ on which $f^n$ is a diffeomorphism in the interior of $J$.
\end{definition}
An endpoint of $J$ is either $0$, $1$ or a preimage of $c$.

For any $x \in [0,1] \setminus \{c\}$ such that $f^n(x) \ne c$ for all $n \in \field{N}$, define the itinerary $\omega(x) \in \{0,1\}^{\field{N}}$ of $x$ as the sequence $\{\omega^0(x), \omega^1(x),\ldots\}$,  such that
\begin{equation}\label{omegas}
\omega^i = \left\{0, \quad f^i(x)<c, \atop 1, \quad f^i(x)>c.\right.
\end{equation}

If one imposes the usual  order $0<1$, then for any two $\omega$ and $\tilde{\omega}$ in $\{0,1\}^{\field{N}}$,  we say that $\omega<\tilde{\omega}$ iff there exists $r \ge 0$ such that $\omega^i = \tilde{\omega}^i$ for all $i<r$ and $\omega^r<\tilde{\omega}^r$.

The limits 
$$\omega(x^+) \equiv \lim_{y \downarrow x}  \omega(y), \quad \omega(x^-)\equiv \lim_{y \uparrow x}  \omega(y)$$
where  $y$'s run through the points which are not the preimages of $c$, exists for all $x \in [0,1]$.

{\it The kneading invariant} $K(f)$ of $f$ is the pair $(K^-(f),K^+(f))=(\omega(c^-),\omega(c^+))$. Hubbard and Sparrow have found in \cite{HS} a condition on the kneading invariant of topologically expansive Lorenz maps. Kneading invariants for a general Lorenz map, not necessarily expansive,  satisfy the following  condition:
$$K^-_0=0, \quad K^+_0=1, \quad \sigma(K^+) \le \sigma^n(K^\pm) \le  \sigma(K^-), \quad n \in \field{N},$$
here $\sigma$ is the shift in $\{0,1\}^{\field{N}}$. Conversely, any sequence as above is a kneading sequence for some Lorenz  map.

A Lorenz map has two critical values 
$$c_1^-=\lim_{x \uparrow c} f(x), \quad c_1^+=\lim_{x \downarrow c} f(x).$$
We will use the notation $c_1^\pm(f)$ whenever we want to emphasize that that a critical value corresponds to a function $f$.

A Lorenz map $f$ with $c_1^+ < c < c_1^-$ is called nontrivial, otherwise $f$ has a globally attracting fixed point. In general, $c^\pm_k$ will denote points in the orbit of the critical values:
$$c^\pm_i=f^{i-1}(c^\pm_1), \quad i \ge 1.$$

\begin{figure}

\begin{center}
\begin{tabular}{c c}
 {\includegraphics[height=53mm,width=50mm,angle=-90]{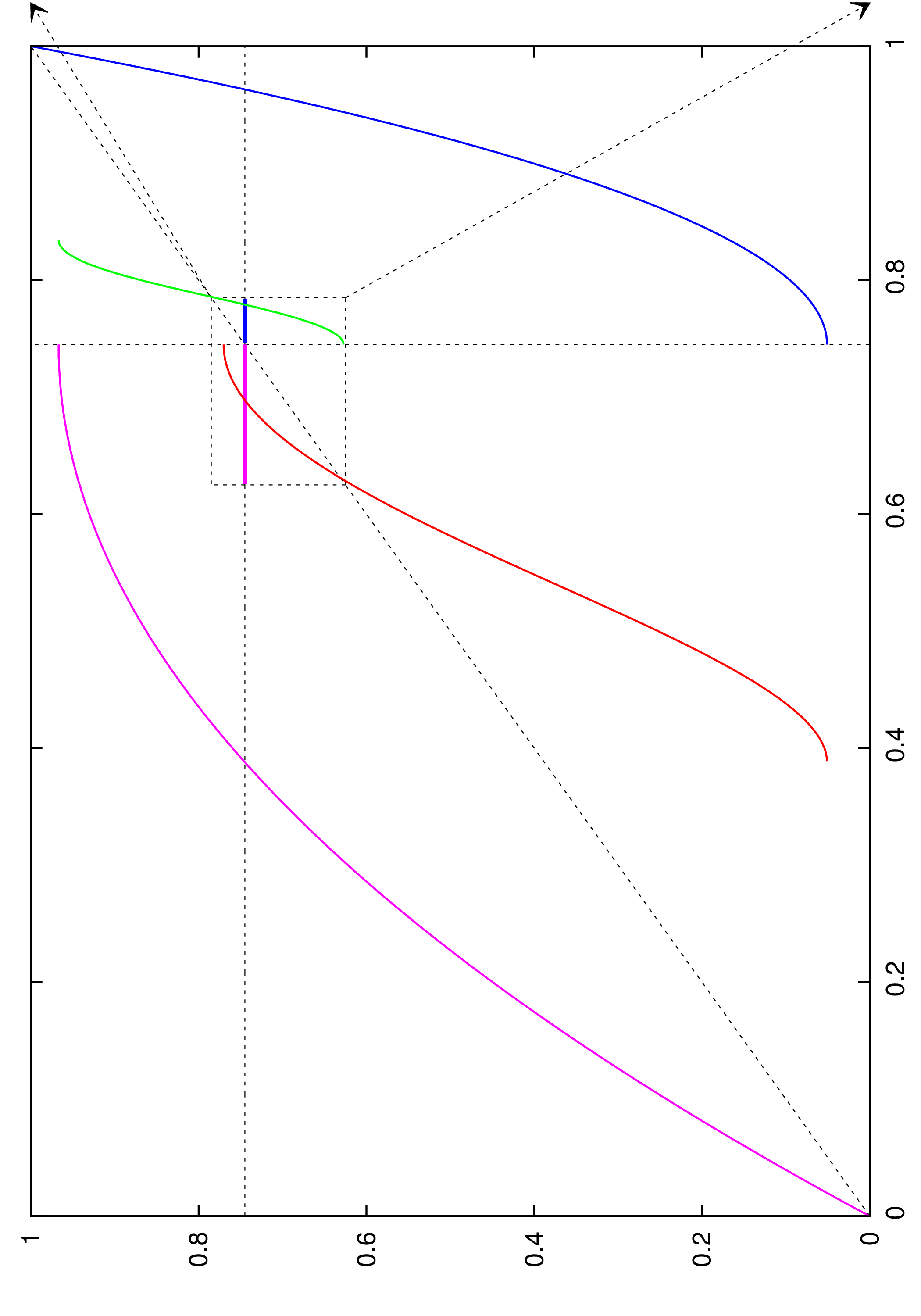}} & {\includegraphics[height=53mm,width=50mm,angle=-90]{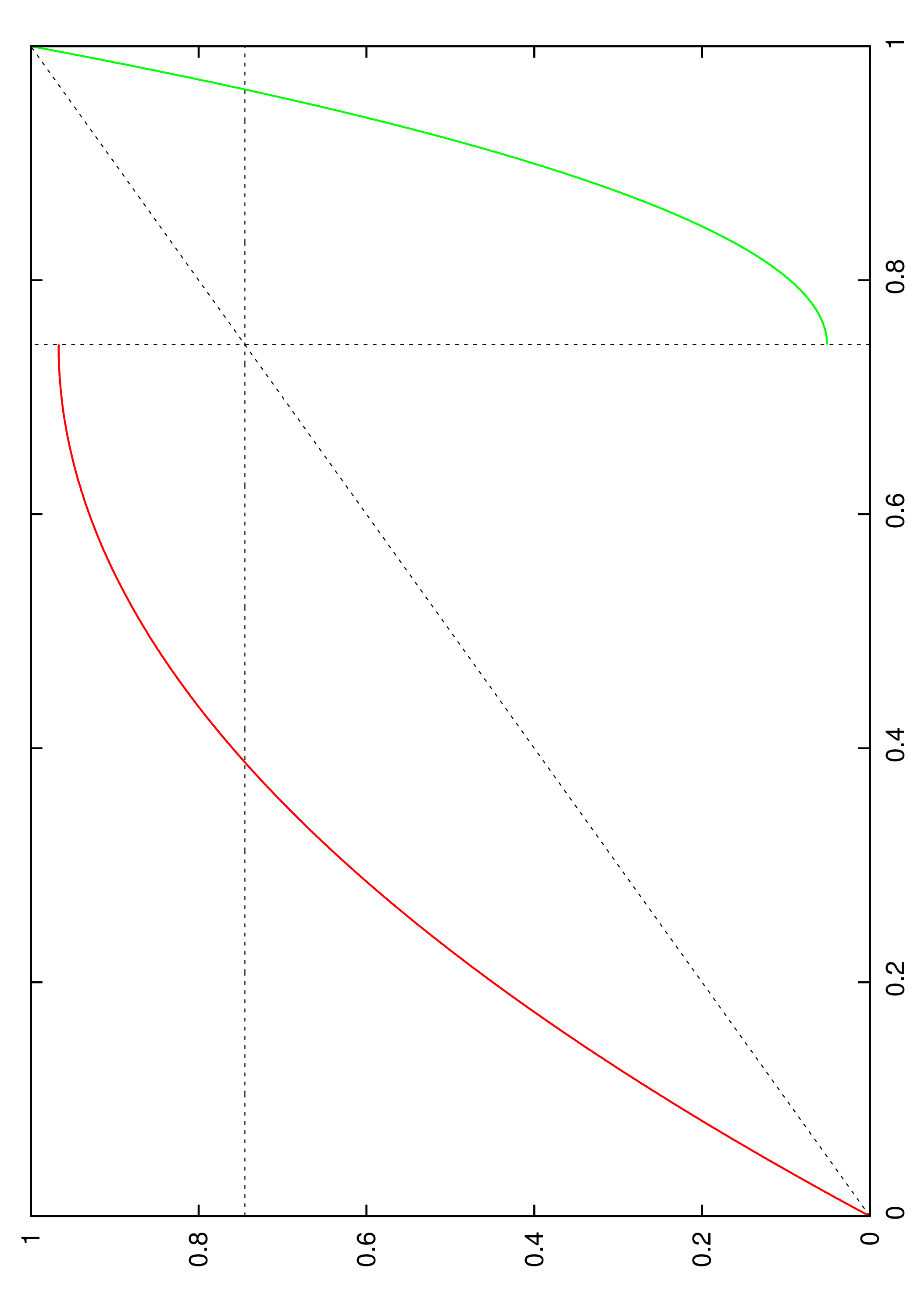}} \\
a) & b)
\end{tabular} 
\caption{a) A Lorenz map $f$ of renormalization type $(01,1000)$ with the critical exponent $\rho=2$; b) $\cR[f]$ } 
\end{center}
\end{figure}

\begin{definition}
A Lorenz map $f$ is called {\it renormalizable} if there exist $p$ and $q$, $0<p<c<q<1$, such that the first return map $(f^n,f^m)$, $n>1, m>1$, of $C=[p,q]$ is affinely conjugate to a nontrivial Lorenz map. Choose $C$ such that it is maximal. The rescaled first return map of such $C \setminus \{c\}$ is called the renormalization of $f$ and denoted $\cR[f]$.
\end{definition}

We will denote 
$$L=[p,c), \quad R=(c,q],$$ 
while the first return map will be denoted $\cP[f]$ and referred to as the prerenormalization. If $f$ is renormalizable, then there exist minimal positive integers $n$ and $m$ such that 
$$\cP[f](x)=\left\{f^{n+1}(x),  \quad x \in L, \atop  f^{m+1}(x),  \quad x \in R,  \right.$$

Then, explicitly,
\begin{equation}\label{ren_op}
\cR[f]=A^{-1} \circ \cP[f] \circ A,
\end{equation}
where $A$ is the affine orientation preserving rescaling of $[0,1]$ onto $C$. We will also use the notation $\tilde{f}$ for the renormalization of $f$.

The intervals $f^i(L)$, $1 \le i \le n$, are pairwise disjoint, and disjoint from $C$. So are the intervals,  $f^i(R)$, $1 \le i \le m$. Since these intervals do not contain $c$, we can associate a finite sequence of $0$ and $1$ to each of these two sequences of intervals:
$$\omega^-\!=\!\{K^-_0,\ldots,K^-_{n}  \},  \  \omega^+\!=\!\{K^+_0,\ldots,K^+_{m}\}, \ \omega=(\omega^-,\omega^+) \in \{0,1\}^{n+1} \times \{0,1\}^{m+1},$$ 
which will be called the {type of renormalization}. The subset of Lorenz maps as in the Definition $(\ref{Lorenz_family})$ which are renormalizable of type $(\omega^-,\omega^+)$ is referred to as the domain of renormalization $\cD_\omega$ (cf.~\cite{MM}).

Let 
\begin{equation}\label{bar_omega}
\bar{\omega}=(\omega_0,\omega_1,\ldots) \in \prod_{i \in \field{N}} \bigotimes \left( \{0,1\}^{n_i+1} \times \{0,1\}^{m_i+1}\right).
\end{equation} 
If $\cR^i[f]$ is $\omega_i$-renormalizable for all $i \in \field{N}$, then $f$ is called infinitely renormalizable of combinatorial type $\bar{\omega}$. The  set of $\omega$-renormalizable maps will be denoted by $\cL_\omega$, and the set of maps $f$ such  that $\cR^i[f]$ is $\omega_i$-renormalizable will be called $\cL_{\bar{\omega}}$, $\bar{\omega}=(\omega_0,\omega_1, ..., \omega_n)$, with $n$ finite or infinite.  If $\bar{\omega}$  is such that  $|\omega_i^\pm|<B$, $i=0,1, \ldots$, for some $0<B<\infty$,  we say that $\bar{\omega}$ is of {\it bounded type}.

We would like to draw the attention of the reader to the position of the indices in our notation: $\omega_i \in \{0,1\}^{\field{N}} \times  \{0,1\}^{\field{N}}$ is a pair of two words, while $\omega^i$ is an integer $0$ or $1$ in a single word (cf. $(\ref{omegas})$).

The combinatorics 
\begin{equation}\label{monotone}
\omega=(0 \overbrace{1 \ldots 1}^{n}, 1 \overbrace{0 \ldots 0}^{m})
\end{equation}
will be called {\it monotone}. The set of all monotone combinatorial types will be denoted $\cM$, while $\cL_\cM$ will denote all Lorenz maps which are $\omega$-renormalizable with $\omega \in \cM$.

\begin{figure}

\begin{center}
 {\includegraphics[height=75mm,width=70mm,angle=-90]{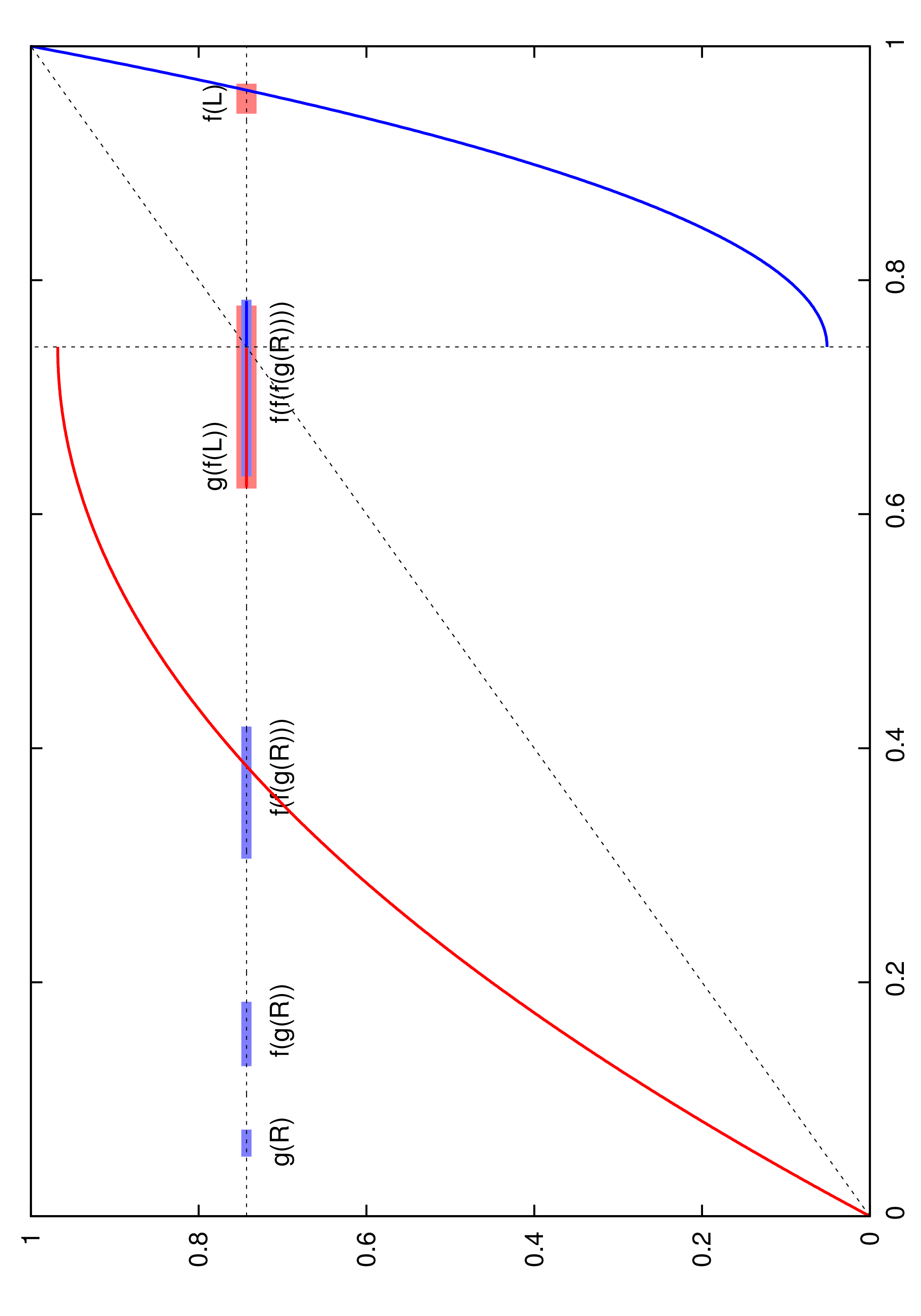}}
\caption{Monotone combinatorics  $(01,1000)$ for a map with the critical exponent $\rho=2$. The two halves of the central interval are given in red and blue, their images under the map in semi-transparent red and blue.}
\end{center}
\end{figure}

Given an integer  $N > 1$, the subset $\cM$  given by all $\omega$'s such that the length of words in $\omega$  satisfies $N \le |\omega^-|$ and  $N \le |\omega^+|$, will be denoted $\cM_N$. Given two integers $M>N>1$, $\cM_{N,M}$ will denote the subset of $\cM$  of all $\omega$'s such that the length of words in $\omega$  satisfies $N \le |\omega^-| \le M$ and  $N \le |\omega^+| \le M$.

Given a subset $\cA \subseteq \cM$,  $\cL_\cA$ will denote all Lorenz maps which are $\omega$-renormalizable with $\omega \in \cA$. We will also use the notation $\cL^S_\cA = \cL^S \cap \cL_\cA$.

\subsection{Statement of results}\label{statement_results}

The main results of our paper are the following proposition and theorems.

\bigskip

\begin{mainprop*}\label{mainprop1} ({\it A priori bounds}).
There exists $\hat{\rho}>2$ such that for every $\rho>\hat{\rho}$  there exists a subset $\cK \subset \cL^1$, relatively compact in $\cL^0$, and a pair of natural numbers $N$ and $M$, such that  $\cR[\cL^S_{\cM_{N,M}} \cap \cK] \subset  \cK$. 
\end{mainprop*}

\bigskip

At this point we were able to prove {\it a priori} bounds only for large $\rho$. The reasons for that will become clear in the proof of the invariance of bounds on the critical point in Proposition $\ref{inv-crpoint}$.  Main Proposition is used to obtain the existence of the periodic points of renormalization:

\bigskip

\begin{thmA*}\label{mainthm1} ({\it Renormalization periodic points}). 
There exists $\hat{\rho}>2$ such that for every $\rho>\hat{\rho}$ and every $\bar{\omega}=(\omega_0, \ldots, \omega_{k-1}) \in \cM_{N,M}^k$, where $N$ and $M$ are  as in in the Main Proposition,  the renormalization operator $(\ref{ren_op})$ has a periodic point in  $\cL^S_{\cM_{N,M}} \cap \cK$ of type $\bar{\omega}$.
\end{thmA*}

\bigskip

The proof of the next Theorem  follows verbatim the one of a similar result in \cite{Win2} (cf Theorem 5.5 in \cite{Win2}), after one establishes  {\it a priori bounds}. We, however, chose to state it as a separate main result since the existence of a Cantor attractor for the dynamics merits a special emphasis.  For completeness, the proof will be included in the Appendix.

\bigskip

\begin{thmB*}\label{mainthm2} ({\it Cantor attractors}).
Let $\rho>\hat{\rho}$, where $\hat{\rho}$ is as in the Main Proposition, and suppose that $\bar{\omega}=(\omega_0, \ldots, \omega_{k-1} \ldots ) \in \cM_{N,M}^\field{N}$, where $N$ and $M$ are  also as in in the Main Proportion.

Consider  $f \in \cL^S_{\bar{\omega}}  \cap \cK$, and  let $\Lambda$ be the closure of the orbits of the critical values. 

Then,
\begin{itemize}
\item[1)] $\Lambda$ is a Cantor set  of a Hausdorff dimension strictly inside $(0,1)$; 
\item[2)] $\Lambda$ is uniquely ergodic;
\item[3)] the complement of the basin  of attraction of $\Lambda$ in $[0,1]$ has zero  Lebesgue measure.
\end{itemize}
\end{thmB*}

\bigskip

The study of renormalizable Lorenz maps was initiated by Tresser et al. (see
e.g. \cite{CCT}).  A more recent work of  Martens and de Melo \cite{MM} produced a series of important results, specifically about the domains of renormalization and the structure of the parameter plane for two-dimensional Lorenz  families.

The work \cite{Win1} presented a computer assisted proof of existence of a renormalization fixed point for the renormalization operator of type  $(\{0,1\},\{1,0,0\})$. The renormalization operator  of this particular type has been later shown  to have  a fixed point in the class of maps analytic on a neighbourhood of the unit interval using only complex analytic techniques in \cite{GaiWin}.

In a more general setting, issues of existence of renormalization periodic points and hyperbolicity have been addressed in \cite{Win2}, where it is proved that the limit set of renormalization, restricted to monotone combinatorics with the return time of one branch being large and much larger than the return time for the other branch, is a Cantor set, and that each point in the limit set has a two-dimensional unstable manifold  (cf Theorem $12.3$ and $12.5$ in \cite{Win2}).  In particular, \cite{Win2} proves equivalents of our Main Proposition (cf Theorem $4.2$ and Theorem $5.1$ in \cite{Win2}) and Theorem A for monotone combinatorial types (cf Theorem $6.1$ in \cite{Win2})  with the following return times:
\begin{equation}\label{wincklers_conds}
[\rho] \le |\omega^-|-1 \le [2 \rho-1], \quad n_- \le |\omega^+|-1 \le n_+,
\end{equation}
where $n_-$ is sufficiently large, and $n_+$ depends on the choice of $n_-$. 

 In comparison, we prove the {\it a priori} bounds for a different class of combinatorial types.  We are able to avoid the disparity of return times evident in $(\ref{wincklers_conds})$, at the price of considering the Lorenz maps with sufficiently flat critical points. In our approach, the exponent $\rho$ has to be not too small (however, we have not computed an explicit lower bound on $\rho$); for each such $\rho$ we demonstrate the existence of two  bounds, upper, $M$, and lower, $N$, such that renormalizable maps of the monotone combinatorics whose length is bounded by $M$ and $N$ admit real {\it a priori} bounds.  The range of the allowed length of the combinatorics is close to 
$$(\ln2 / \ln \rho +1) \rho <  |\omega^-|, |\omega^+| < 2 \rho.$$
In comparison to $(\ref{wincklers_conds})$, the length of combinatorics in our proofs is similar to the one for the ``shorter'' branch in \cite{Win2}. 

Another difference between the class of maps considered in \cite{Win2}, and in our paper, is the location of the critical point. The renormalization invariant set $\cK$ considered in \cite{Win2} consists of maps whose critical point is located close to one of the endpoints of the unit interval.  In our case, the critical point is allowed to be in an interval symmetric with respect to $1/2$, and approaching $[0,1]$ as $\rho$ increases: specifically, $c \in [\rho^{-\beta},1-\rho^{-\beta}]$ where $\beta <1$ is close to $1$.

Our proof uses ideas similar to those in \cite{Win2}, however, because of the differences in the combinatorial types and the set $\cK$, many technical details of the proofs are quite different. This is especially evident in the proof of the invariance of the set $K$, Propositions $\ref{inv_distortion}$, $\ref{inv-crpoint}$ and $\ref{apriori}$.

\bigskip

\section{Preliminaries}
\subsection{The Koebe Principle}

We  start by quoting the Koebe Principle which is of a fundamental importance in real dynamics (see, ex. \cite{dMvS}). We will say that an interval $V$ is a $\tau$-scaled neighbourhood of $U \subset V$, if both components of $V\setminus U$ have length at least $\tau \cdot U$. 

\begin{KPr*} \label{KP}
Let $J \subset T$ be intervals, and $f: T \mapsto f(T)$ be a $C^{3}$-diffeomorphism with $S_f<0$. If $f(T)$ contains a $\tau$-scaled neighbourhood of $f(J)$, then  
$$\left(\tau \over 1+\tau \right)^2 \le {f'(x) \over f'(y) } \le \left(1+\tau \over \tau \right)^2, \quad x,y  \in J.$$
\end{KPr*}

\subsection{Distortion and nonlinearity}
Let $C^k(A;B)$ be the set of $k$-continuously differentiable  maps from $A$ to $B$. We denote  $\cD^k(A;B) \subset C^k(A;B)$  the subset of orientation preserving homeomorphisms whose inverse lies in $C^k(A;B)$.  We will use the notation $\cD^k$ and $C^k$ whenever $A=B=[0,1]$.

\begin{definition}
The {\it nonlinearity} operator $N: \cD^2(A;B) \mapsto C^0(A;\field{R})$ is defined as
$$N_\phi=\left(\log{  \phi' }\right)',$$
while
$$N_{\phi}(x)={\phi''(x) \over \phi'(x)}$$
is the nonlinearity of $\phi$ at point $x$.
\end{definition}

\begin{definition}
Given  $\phi \in \cD^1(A;B)$, the quantity
$${\rm dist}[\phi]=\max_{x,y \in  A}  \ln \left( {\phi'(y) \over \phi'(x) } \right)$$
is called the {\it distortion} of $\phi$.
\end{definition}

Notice, that 
$$\int_x^{y} N_\phi(t) d t= \ln{ \phi'(y) \over \phi'(x)   }.$$

The following Lemma results from a straightforward computation.

\begin{lemma}
The nonlinearity operator $N: \cD^2(A;B) \mapsto C^0(A;\field{R})$ is a bijection. In  the case $A=B=[0,1]$, the inverse is  defined as
\begin{equation}\label{inv_non_op}
N^{-1}_\phi (x)={ \int_0^x \exp\left\{\int_0^r\phi(t) d t  \right\} d r  \over \int_0^1 \exp\left\{\int_0^r\phi(t) d t  \right\} d r  }.
\end{equation}
\end{lemma}

One can turn $\cD^2(A;B)$ into a Banach space using the nonlinearity operator. Specifically, for $\phi$, $\psi$ in $\cD^2(A;B)$ and $a,b \in \field{R}$, the linear structure and the norm are defined via
\begin{eqnarray} 
\label{lin_struct} a \phi+ b \psi&=&N^{-1}_{a  N_{\phi}+b  N_{\psi}},\\
\|\phi\|&=&\sup_{x \in A} \left|N_{\phi}(x) \right|.
\end{eqnarray}

Finally, we give a list of useful bounds on derivatives and distortion in $\cD^2(A;B)$  in terms on the nonlinearity (see \cite{Win2}, Lemma B.10, Lemma B11, or \cite{dMvS} for the proofs).

\begin{lemma}
If $\phi, \psi  \in  \cD^2(A;B)$  then, for all $x,y \in  A$,
\begin{equation}\label{n1}
e^{-|y-x|  \| \phi \| } \le {\phi'(y) \over \phi'(x)} \le e^{|y-x|  \| \phi \| },
\end{equation}
\begin{equation}\label{n2}
{|B| \over |A|}  e^{- \| \phi \| } \le {\phi'(x)} \le {|B| \over |A|}e^{ \| \phi \| },
\end{equation}
\begin{equation}\label{n3}
e^{-\| \phi-\psi \| } \le {\phi'(x) \over \psi'(x)} \le e^{\| \phi-\psi \| }.
\end{equation}
\end{lemma}


We will introduce two subsets of Lorenz maps, defined via conditions on their distortion and  critical points.

\begin{definition}
Given a real constants $\pi>0$, we set
\begin{equation}\label{set1}
\cK^\pi  \equiv  \left\{f \in \cL^1: {\rm dist}[\psi] \le \pi, {\rm dist}[\phi] \le \pi \right\}.
\end{equation}

Given real constants $\pi>0$, $\varepsilon >0$, set
\begin{equation}\label{set2}
 \cK^\pi_{\vareps}  \equiv  \left\{f \in \cK^\pi \subset \cL^1: c(f) \in [\vareps, 1-\vareps] \right\}.
\end{equation}
\end{definition}

The reason for the introduction of these sets is the following compactness result.

\bigskip

\begin{cor}
Given $\pi>0$ and $\vareps>0$, the set $\cK^\pi_\vareps$ is relatively compact in $\cL^0$.
\end{cor}
\begin{proof}
Recall that $\cL^1$ is isomorphic to $[0,1]^2 \times (0,1) \times \cD^1 \times \cD^1$. Since $c$ is bounded away from $0$ and $1$ by a constant, it is, therefore,  contained in a compact  subset  of $(0,1)$. Consider  the set  
$$\cB=\left\{(\phi,\psi) \in \cD^1 \times \cD^1: {\rm dist}[\phi]  \le \pi, {\rm  dist}[\psi]\le\pi \right\}.$$

Any sequence from $\cB$ is equicontinous since $|\phi(y)-\phi(x)|  \le e^\pi |y-x|$, and, clearly, uniformly  bounded, therefore by the Arzel\`a-Ascoli theorem $\cB$ is sequentially compact in $C^0$ topology. The claim follows.
\end{proof}
\medbreak

\subsection{Monotone combinatorics}

We will quote a lemma from \cite{Win2}  (Lemma $2.11$) which gives the formulae for the factors of the renormalization of a  Lorenz map in $\cL_\cM$. Let $I$ be an interval and $g_I$ be an orientation preserving diffeomorphism. We denote the affine transformation that takes $[0,1]$ onto $I$ as $\xi_I$. Define the {\it zoom operator}:
\begin{equation}\label{zoom}
Z(g;I)=\xi^{-1}_{g(I)} \circ g \circ \xi_I.
\end{equation}

\begin{lemma}\label{renorm_parts}
If $f=(u,v,c,\phi,\psi)$ is renormalizable of monotone combinatorics, then
$$\cR[f]=(\tilde{u},\tilde{v},\tilde{c},\tilde{\phi},\tilde{\psi})$$
is given by
\begin{equation}
\label{uvc'} \tilde{u}={|Q(L)| \over |U|}, \quad \tilde{v}={|Q(L)| \over |V|}, \quad \tilde{c}={|L| \over |C|}, 
\end{equation} 
\begin{equation}
 \label{phipsi'} \tilde{\phi}=Z(\bar{\phi}; U), \quad  \tilde{\psi}=Z(\bar{\psi}; V), \quad \bar{\phi}=f_1^n \circ \phi, \quad  \bar{\psi}=f_0^m \circ \psi,
\end{equation} 
where $U=\phi^{-1} \circ f_1^{-n}(C)$, $V=\psi^{-1} \circ f_0^{-m}(C)$.
\end{lemma}

\bigskip

\section{Estimates for Lorenz maps with monotone combinatorics}

In this Section we will obtain bounds on the critical points, critical values and lengths of the central subintervals $L$ and $R$ for Lorenz maps with monotone combinatorics whose diffeomorphic coefficients have bounded distortion.

The main result of this section is Lemma $\ref{lem3}$ which gives upper and lower bounds on lengths of renormalization intervals $L$ and $R$ in terms of the exponent $\rho$, the upper bound on the distortion $\pi$, the critical point $c$, the critical values $c^\pm_1$ and the lengths of the combinatorics $\omega_\pm$. The rest of the lemmas in this section are preparatory estimates needed for the proof of Lemma $\ref{lem3}$. We will explain the ideas involved in proving Lemma $\ref{lem3}$ immediately before the Lemma.

Denote
\begin{eqnarray}
\nonumber f_1^{\circ n} \circ f_0(L) &\equiv& I \equiv L_{n+1} \equiv [p,c^-_{n+1}), \quad  f_0(L)=L_1,  \quad  f_1^{\circ k} \circ f_0(L)=L_{k+1}, \\
 \nonumber f_0^{\circ m} \circ f_1(R) &\equiv& J \equiv R_{m+1} \equiv(c^+_{m+1},q], \quad  f_1(R)=R_1,  \quad  f_0^{\circ k} \circ f_1(R)=R_{k+1}.
\end{eqnarray}

We will mention the following simple lemma.

\begin{lemma}  \label{psiphi_ext_lemma}
Suppose $f \in \cL^1$, then $(f_0^{-1})^{\circ n}$  and $(f_1^{-1})^{\circ m}$  are diffeomorphisms on  $(0,c^-_1)$ and $(c^+_1,1)$.
\end{lemma}
\begin{proof}
The branches  $f_0^{-1}$  and $f_1^{-1}$  map the intervals  $(0,c^-_1)$ and $(c^+_1,1)$, respectively, diffeomorphically into themselves.
\end{proof}

We will continue with a sequence of lemmas which will prepare us for the construction of {\it a priori} bounds ---  construction of a relatively compact set invariant under renormalization.

First of all, we will need  simple bounds on the difference of $f_0$ and $f_1$ at two points of the domain.

\begin{lemma}\label{anbn}
Suppose that ${\rm dist}[\phi] \le \pi$,  ${\rm dist}[\psi] \le \pi$, then 
\begin{equation} \label{anbn1}
 {e^{-\pi} \rho c_1^- \over c}  (x-y) \left({c-x \over c }\right)^{\rho-1} \le f_0(x) - f_0(y) \le {e^{\pi} \rho c_1^- \over c}  (x-y) \left({c-y \over c }\right)^{\rho-1},
\end{equation}
for any $x>y$ in $[0,c)$, and
\begin{equation} \label{anbn2}
{e^{-\pi} \rho (1-c_1^+) \over \mu}  (x-y) \left({y-c \over \mu }\right)^{\rho-1} \le f_1(x) - f_1(y) \le {e^{\pi} \rho (1-c_1^+) \over \mu}  (x-y) \left({x-c \over \mu }\right)^{\rho-1},
\end{equation}
for any $x>y$ in $(c,1]$.

\end{lemma}
\begin{proof}
Notice that the average derivative of $\phi$ on $(0,u)$ is $c_1^+/u$, therefore, the derivative $\phi'(x)$ at any point in $(0,u)$ is bounded as 
\begin{equation}\label{phider}
{c_1^- \over u} e^{-\pi} \le \phi'(x) \le {c_1^- \over u} e^{\pi}.
\end{equation}
Similarly, for $x \in (1-v,1)$.
\begin{equation}\label{psider}
{1-c_1^+ \over v} e^{-\pi} \le \psi'(x) \le {1-c_1^+ \over v} e^{\pi}.
\end{equation}
Therefore, we get for $x>y$ in $[0,c)$
$$f_0(x) - f_0(y) \le {c_1^- \over u} e^{\pi} \rho {u \over c} (x-y) \left({c-y \over c }\right)^{\rho-1} = {e^{\pi} \rho c_1^- \over c}  (x-y) \left({c-y \over c }\right)^{\rho-1}.$$
The lower bound is obtained as follows: 
$$f_0(x) - f_0(y) \ge {c_1^- \over u} e^{-\pi} \rho {u \over c} (x-y) \left({c-x \over c }\right)^{\rho-1} = {e^{-\pi} \rho c_1^- \over c}  (x-y) \left({c-x \over c }\right)^{\rho-1}.$$

Bounds on the difference of $f_1$ can be obtained in a similar way.
\end{proof}
\medbreak

For the sake of brevity, let us introduce the following notation:
\begin{equation}
\label{alpha} \alpha \equiv {e^{-\pi}\over \rho}, \quad \eta \equiv {e^{-\pi} \mu \over (1-c_1^+) \rho}, \quad \kappa \equiv {e^{-\pi} c \over c_1^- \rho}, \quad \gamma \equiv {e^{2 \pi} \over \rho},  \quad \nu \equiv {\mu \over (1-c_1^+)^{1 \over \rho}}, \quad \xi \equiv  {c \over (c_1^-)^{1 \over \rho}}.
\end{equation}

Since $R \subset f^{m+1}(R)$, we have that $f_0^{-1}(c) \in f^m(R)$, and, therefore,  for monotone combinatorics $f_0^{-1}(c)>c^+_1$. Similarly, $f_1^{-1}(c)<c_1^-$.  The next lemma uses this fact, and provides a lower bound on the length of the intervals $[f^{-1}_0(c),p]$ and $[q,f^{-1}_1(c)]$, which is also a lower bound on the length of the intervals $[c^+_1,p]$ and $[q,c^-_1]$.

\begin{lemma}\label{delt}
Let $f \in \cK^\pi \cap \cL_\omega$ for some $\pi >0$ and $\omega=(\omega_-,\omega_+) \in \cM$ with  $|\omega_-|=n+1$, $|\omega_+|=m+1$. Then
\begin{eqnarray}
\label{pf0} |p-f_0^{-1}(c)| &\ge&   \left(\kappa {\left({c \over c-c_1^+ } \right)^{\rho-1}}  \left( \nu^{\rho \over \rho-1}  e^{-\pi  \over \rho-1} \right) \right)^{\rho^n \over \rho^n-1} \equiv \Delta,\\
\label{qf1} |q-f_1^{-1}(c)| &\ge&    \left(\eta {\left({\mu \over c_1^--c } \right)^{\rho-1}}  \left( \xi^{\rho \over \rho-1}  e^{-\pi  \over \rho-1} \right) \right)^{\rho^m \over \rho^m-1} \equiv \Theta.
\end{eqnarray}
\end{lemma} 
\begin{proof}
We will first demonstrate that
\begin{equation}\label{induc}
f_1^{-n}(x) \ge c+\nu^{\rho \over \rho-1}  e^{-{\pi  \over \rho-1}} \left(x-c_1^+  \right)^{1 \over \rho^n}
\end{equation}
for all  $x > c_1^+$.  To prove $(\ref{induc})$ we use the following expressions for the inverse branches of a Lorenz map:
\begin{eqnarray}
\label{f0m1} f_0^{-1}(x)&=&c-c \left({|\phi^{-1}([x,c_1^-])| \over  |\phi^{-1}([0,c_1^-])|  }\right)^{1 \over \rho}= c-c \left({u-\phi^{-1}(x) \over  u }\right)^{1 \over \rho}, \\ 
\label{f1m1} f_1^{-1}(x)&=&c+\mu \left(1-{|\psi^{-1}([x,1])| \over  |\psi^{-1}([c_1^+,1])|}  \right)^{1 \over \rho}= c+\mu \left(1-{1-\psi^{-1}(x) \over  v }\right)^{1 \over \rho}. 
\end{eqnarray}

Start with 
$$f_1^{-1}(x) \ge c+\mu \left( e^{-\pi} {x-c_1^+ \over 1-c_1^+}   \right)^{1 \over \rho}= c+\mu \left( {e^{-\pi} \over 1-c_1^+} \right)^{1 \over \rho} \left(x-c_1^+   \right)^{1 \over \rho},$$
for $x > c_1^+$, and use induction on this inequality to obtain
$$f_1^{-n}(x) \ge c+ \left(\mu  \left({e^{-\pi} \over 1-c_1^+} \right)^{1\over \rho} \right)^{1+\ldots \rho^{-(n-1)}} \left( x-c_1^+\right)^{1 \over \rho^n} \ge c+\nu^{\rho \over \rho-1}  e^{-\pi  \over \rho-1} \left(x-c_1^+   \right)^{1 \over \rho^n}.$$

According to Lemma $\ref{anbn}$:
$$|f_0(p)-c| \le  \kappa^{-1} |p-f_0^{-1}(c)|  \left({c-f_0^{-1}(c) \over c} \right)^{\rho-1}.$$

On the other hand, $f_0(p)=f_1^{-n}(p)$, and according to $(\ref{induc})$,
\begin{equation}
f_1^{-n}(p) \ge c+\nu^{\rho \over \rho-1}  e^{-\pi  \over \rho-1} \left( p-c_1^+   \right)^{1 \over \rho^n}.
\end{equation}
Therefore,
$$f_0(p)-c = f_1^{-n}(p)-c \ge \nu^{\rho \over \rho-1}  e^{-\pi  \over \rho-1} | p-c_1^+|^{1 \over \rho^n},$$
and
\begin{equation}
\label{ineq1}\kappa^{-1}|p-f_0^{-1}(c)|  {\left({c-f_0^{-1}(c) \over c} \right)^{\rho-1}} \ge \nu^{\rho \over \rho-1}  e^{-\pi  \over \rho-1} | p-c_1^+|^{1 \over \rho^n} \ge  \nu^{\rho \over \rho-1}  e^{-\pi  \over \rho-1} | p-f_0^{-1}(c)|^{1 \over \rho^n},
\end{equation}
which results in the required bound $(\ref{pf0})$.

The bound on  $|q-f_1^{-1}(c)|$ is obtained in a similar way.
\end{proof}
\medbreak

Lower bounds on the differences  $|p-f_0^{-1}(c)|$ and  $|f_1^{-1}(c)-q|$ can be used to bound $c_1^-$ and $1-c_1^+$ from below.

\begin{lemma}\label{lem2}
Let $f \in \cK^\pi \cap \cL_\omega$ for some $\pi >0$ and $\omega=(\omega_-,\omega_+) \in \cM$ with  $|\omega_-|=n+1$, $|\omega_+|=m+1$. Then,
\begin{equation}\label{c1pm}
c_1^+ \ge  {\kappa^m  \Delta \over 1-\kappa^m}, \quad 1-c_1^- \ge  {\eta^n  \Theta \over 1-\eta^n}.
\end{equation}
\end{lemma}
\begin{proof}
To get the lower bound on $c_1^+$ we notice that the derivatives of the inverse branches of $Q(x)$ (formulae $(\ref{f0m1})$ and $(\ref{f1m1})$ with $\phi=\psi={\rm id}$) are increasing functions, while the derivatives of $\phi$ and $\psi$ are bounded as in $(\ref{phider})$ and $(\ref{psider})$. This can be used to get a straightforward bound  
$$(f^{-1}_0)'(x) \ge  {e^{-\pi} c \over c_1^- \rho}=\kappa,$$
for all $0<x<c_1^-$. Therefore,
$$f^{-m}(p) \ge \left( Df^{-1}(0) \right)^m p \ge \kappa^m p,$$
so 
$$p \ge c_1^+ + \Delta  \ge \kappa^m p +\Delta  \implies p  \ge { \Delta \over 1-\kappa^m},$$
and
$$c_1^+ \ge  {\kappa^m  \Delta \over 1-\kappa^m}.$$

The lower bound on $1-c_1^-$ is obtained in a similar way.
\end{proof}
\medbreak

We will now turn our attention to the bounds  on $L$ and $R$. 

The first key idea in obtaining an upper bound, for example, on the size of the interval $L$, is that the orbit of this interval, $f^i(L)$, $i=1,..n$, is contained in  $(f_1^{-1}(c),c_1^-)$, therefore, the length $|(f_1^{-1}(c),c_1^-)|=c_1^--f_1^{-1}(c)$ bounds the sum $\sum_{i=1}^n|f^i(L)|$ from above. The second observation has to do with the form of the map $f_1$ along the orbit of $L$: to bound the lengths of intervals $f^i(L)$ in terms of $|L|$, one uses first that $f\arrowvert_L$ is essentially a power map, while $f\arrowvert_{f^i(L)}$, $i=1, \ldots,n$, is essentially linear whose derivative can be bounded in terms of $\rho$, $c$, $c_1^+$, $n$ and $\pi$. These two observations put together produce an upper bound on $|L|$.

To obtain lower bounds, say, again, on $|L|$, we first notice that $L$ covers itself under the return map $f^{n+1}$, therefore $|L| \le |f^{n+1}(L)|$. Next, one estimates the change of the size of $f(L)$ under $f^n$, using bounds from Lemma $\ref{anbn}$, and, again, the fact that $|f(L)| \sim |L|^\rho$ to obtain that $|L| \le |f^{n+1}(L)| \le C |L|^\rho \implies |L| \ge C^{-{1 \over \rho-1}}$, where $C$ is some constant that depends on $c$, $c_1^\pm$, $n$ and $\pi$.

\begin{lemma}\label{lem3}
Let $f \in \cK^\pi \cap \cL_\omega$ where $0<2 \pi<\ln \rho$ and $\omega=(\omega_-,\omega_+) \in \cM$ with  $|\omega_-|=n+1$, $|\omega_+|=m+1$. Then there exist a constant $K$, such that 
\begin{eqnarray}
 \label{upperL1} |L| &\le&   \left((c_1^--q)  {c^\rho e^{\pi}  \over c_1^-} \right)^{1\over \rho+1}  \left(   { \gamma^{-1}-1 \over \gamma^{-n}-1   } \right)^{1 \over  \rho+1}, \\
\label{upperL2}  
|L|  &\le& {\mu^{\rho \over \rho+1} e^{2 \pi {2 \rho-1 \over \rho (\rho+1)} } c^{\rho-1 \over \rho+1} \over \rho^{1 \over \rho+1}   }  {1 \over |R|^{\rho-1 \over \rho}}   \left(\gamma^{-1} -1 \over \gamma^{-m}-1   \right)^{\rho-1 \over \rho (\rho+1)} \left( { 1-\gamma \over \gamma^{2-n}-\gamma  } \right)^{1 \over \rho+1} \left( 1 \over 1+{e^{-\pi} |L| \over |R|^{\rho+1} \left({ \gamma^{-m}-1 \over \gamma^{-1}-1}   \right)}  \right)^{\rho-1 \over \rho (\rho+1)}, \\
\label{lowerL}   |L| &\ge& \left( {e^{-\pi} c^\rho \over  c_1^-}  \eta^n \right)^{1 \over \rho-1}  \exp \left( K {\eta^n \Theta \over \mu( 1-\eta^n)} \sum_{k=1}^n \left({e^{-2 \pi} \over \eta} {(\Theta+|R|)^{\rho-1} \over \mu^{\rho-1}}  \right)^{k-1} \right),
\end{eqnarray}
and
\begin{eqnarray}
\label{upperR1}  |R| &\le&  \left((p-c_1^+)  {\mu^\rho e^{\pi}  \over (1-c_1^+)} \right)^{1\over \rho+1} \left(   { \gamma^{-1}-1 \over \gamma^{-m}-1   } \right)^{1 \over  \rho+1}, \\
\label{upperR2}  |R| &\le& {c^{\rho \over \rho+1} e^{2 \pi {2 \rho-1 \over \rho (\rho+1)} } \mu^{\rho-1 \over \rho+1} \over \rho^{1 \over \rho+1}   }  {1 \over |L|^{\rho-1 \over \rho}}    \left(\gamma^{-1} -1 \over \gamma^{-n}-1     \right)^{\rho-1 \over \rho (\rho+1)} \left( { 1-\gamma \over \gamma^{2-m}-\gamma  } \right)^{1 \over \rho+1} \left( 1 \over 1+{e^{-\pi} |R| \over |L|^{\rho+1} \left({ \gamma^{-n}-1 \over \gamma^{-1}-1}   \right)}  \right)^{\rho-1 \over \rho (\rho+1)}, \\ 
\label{lowerR}  |R| &\ge& \left( {e^{-\pi} \mu^\rho \over  1-c_1^+}  \kappa^m \right)^{1 \over \rho-1}  \exp \left( K {\kappa^m  \Delta \over c (1-\kappa^m) } \sum_{k=1}^m \left({e^{-2 \pi} \over \kappa} {(\Delta+|L|)^{\rho-1} \over c^{\rho-1}}  \right)^{k-1}  \right).
\end{eqnarray}
\end{lemma}

\begin{proof}
\noindent {\it 1)  Upper bounds}.  Denote $p_i=f^i(p)$ and $q_i=f^i(q)$ (notice, $p_{n+1}=p$ and $q_{m+1}=q$), and, as before, $c_i^\pm=f^{i-1}(c_1^\pm)$.  Suppose, point $x_1$ is in the interval $L_1$, and denote points in the orbit of $x_1$ as $x_k$: $x_k=f_1^{k-1}(x_1)$. Then, according to $(\ref{induc})$,
\begin{equation}
\label{eq:xk}x_{k} \equiv f_1^{-(n-k)}(x_n)  \ge c+\nu^{\rho \over  \rho-1} e^{-{\pi \over \rho-1}} \left( x_n -c_1^+ \right)^{1 \over \rho^{n-k}}\equiv \tilde{x}_k,
\end{equation}
and one gets for all $n+1>k>0$
\begin{eqnarray}
\nonumber (f_1^{-1})'(x_{k+1}) \!\!\! &\le&   \!\!\!  {\mu e^{\pi} \over \rho (1-c_1^+)} \left( e^{-\pi} {x_{k+1}-c_1^+ \over 1-c_1^+} \right)^{1-\rho \over \rho}  \!\!\!\! \le  {\mu e^{\pi} \over \rho (1-c_1^+)}  (1-c_1^+)^{\rho-1 \over \rho} \!\left( e^{-\pi} (x_{k+1}-c )\right)^{1-\rho \over \rho} \\
\nonumber  &\le& {\mu e^{\pi} \over \rho}  (1-c_1^+)^{-{1 \over \rho}} \left(\! e^{-\pi} \! \left( {\mu \over (1-c_1^+)^{1 \over \rho}  } \right)^{\rho \over \rho-1} \!\! e^{-{\pi \over \rho-1}} (x_n-c_1^+)^{1 \over \rho^{n-k-1}}  \! \right)^{1-\rho \over \rho} \\
\nonumber   \!\!\! & \le &   \!\!\! {e^{2 \pi}  \over \rho} \left(x_n -c_1^+  \right)^{1-\rho \over \rho^{n-k}}=\gamma \left(x_n -c_1^+  \right)^{1-\rho \over \rho^{n-k}},      
\end{eqnarray}
where we have used $(\ref{eq:xk})$ going from the first line to the second. Therefore,
$$f_1'(x_{k})=  \left((f_1^{-1})'(x_{k+1}) \right)^{-1} \ge  \gamma^{-1} \left({x_n-c_1^+} \right)^{\rho-1 \over \rho^{n-k}}.$$
We can now see that
\begin{equation} 
\label{eq:Df1} f_1'(x_k) \ge  \gamma^{-1} \left({p_n-c_1^+} \right)^{\rho-1 \over \rho^{n-k}}
\end{equation} 
for all $x_k \in L_k$. Therefore,
\begin{eqnarray}
\nonumber L_k&=&|c_{k}^--p_{k}|  \ge |p_1-c_1^-| \prod_{i=1}^{k-1} \min_{x \in L_i} f_1'(x) \ge |p_1-c_1^-| \prod_{i=1}^{k-1} \gamma^{-1} \left({p_n -c_1^+} \right)^{\rho-1 \over \rho^{n-i}} \\
 \nonumber &\ge& |p_1-c_1^-| \gamma^{1-k}  \left({p_n-c_1^+} \right)^{1\over \rho^{n-k}}.
\end{eqnarray}
Notice, that for monotone combinatorics all intervals $L_k=f^k(L)$, $1 \le k \le n-1$, are contained in the interval $(f_1^{-1}(c),c_1^-)$, while the intervals $R_k=f^k(R)$, $1 \le k \le m-1$, are all contained in $(c_1^+,f_0^{-1}(c))$.  Therefore,
\begin{equation}\label{c1minc}
c_1^--f_1^{-1}(c)  > \sum_{k=1}^{n-1}|c_{k}^--p_{k}| \ge |p_1-c_1^{-}| \sum_{k=1}^{n-1}   \gamma^{1-k}  \left({p_n-c_1^+} \right)^{1\over \rho^{n-k}} \ge {c_1^- \over e^{\pi}}  \left({|L|  \over c}  \right)^\rho  \sum_{k=1}^{n-1}  \gamma^{1-k}  |L|^{1\over \rho^{n-k}},
\end{equation}
where we have used that $L \subset (c_1^+,p_n)$ in the last inequality. We can now use the fact that $\gamma^{-1}=\rho / e^{2 \pi} >1$ for all $\pi$ as in the hypothesis of the Lemma, to simplify the above expression. 
\begin{equation}
\label{bound_int}   c_1^--f_1^{-1}(c) \ge {c_1^- \over e^{\pi}}  \left({|L|  \over c}  \right)^\rho    \gamma^{2-n} |L| \sum_{k=0}^{n-2}  \gamma^k={c_1^- \over e^{\pi} c^\rho}  |L|^{\rho+1}   { \gamma^{2-n}-\gamma \over 1-\gamma   }.
\end{equation}
Similarly, since all intervals $L_k$,  $1 \le k \le n$, are contained in the interval $(q,c_1^-)$, while the intervals $R_k$, $1 \le k \le m$, are all contained in $(c_1^+,p)$, we also have  
\begin{eqnarray}
\nonumber c_1^--q &\ge& {c_1^- \over e^{\pi}}  \left({|L|  \over c}  \right)^\rho    \gamma^{1-n} |L| \sum_{k=0}^{n-1}  \gamma^k={c_1^- \over e^{\pi} c^\rho}  |L|^{\rho+1}   { \gamma^{1-n}-\gamma \over 1-\gamma},
\end{eqnarray}
and the upper bound $(\ref{upperL1})$ from the claim follows. The bound $(\ref{upperR1})$ on $R$ is obtained in a similar way.

To obtain bound $(\ref{upperL2})$, we will return to $(\ref{c1minc})$ and  find an upper bound on $c_1⁻-f_1^{-1}(c)$, using  $(\ref{upperL1})$. Notice, that by Lemma $\ref{anbn}$  and by $(\ref{upperR1})$,
\begin{eqnarray}
\nonumber \mu&>&|f_1(c_1^-)-c|=|f_1(c_1^-)-f_1(f_1^{-1}(c))| \\
\nonumber &\ge& {e^{-\pi} \rho (1-c_1^+) \over \mu   } (c_1^--f_1^{-1}(c)) \left({f_1^{-1}(c)-c  \over \mu  }\right)^{\rho-1} \\
\nonumber &\ge & {e^{-\pi} \rho (1-c_1^+) \over \mu   } (c_1^--f_1^{-1}(c)) \left(e^{-\pi} {c-c_1^+ \over 1-c_1^+} \right)^{\rho-1 \over \rho}  \\
\nonumber &\ge& {e^{-\pi} \rho (1-c_1^+) \over \mu   } (c_1^--f_1^{-1}(c)) \left(e^{-\pi} {|R|^{\rho+1}  \left(\gamma^{-m}-1 \over \gamma^{-1} -1 \right)  \over \mu^\rho e^\pi  } +e^{-\pi} {|L| \over 1-c_1^+}\right)^{\rho-1 \over \rho}   \implies \\
\nonumber c_1^--f_1^{-1}(c) & \le & {\mu e^\pi \over \rho} \left(\mu^\rho e^{2 \pi} \over |R|^{\rho+1}  \left(\gamma^{-m}-1 \over \gamma^{-1} -1 \right) +{\mu^\rho e^\pi |L| \over 1-c_1^+}   \right)^{\rho-1 \over \rho}.
\end{eqnarray}
We use the above inequality together with $(\ref{bound_int})$,
\begin{equation}
{\mu e^\pi \over \rho} \left(\mu^\rho e^{2 \pi}  \over |R|^{\rho+1}  \left(\gamma^{-m}-1 \over \gamma^{-1} -1 \right) + \mu^\rho e^\pi |L|   \right)^{\rho-1 \over \rho} \ge {c_1^- \over e^{\pi} c^\rho}  |L|^{\rho+1}   { \gamma^{2-n}-\gamma \over 1-\gamma }
\end{equation}
which results in bound $(\ref{upperL2})$. The bound $(\ref{upperR2})$ is obtained in a similar way.

\bigskip

\noindent {\it 2) Lower bounds}.  We will use the fact that  $L \subset L_{n+1}=f^{n+1}(L)$, or  
$$ |L| \le |f^n_1(p_1)-f_1^n(c_1^-)|.$$
Then, according to  Lemma $\ref{anbn}$,
\begin{eqnarray}
\nonumber |L| &\le& |f_1(p_n)-f_1(c_n^-)|\le  {e^\pi (1-c_1^+) \rho \over \mu } \left|p_n-c_n^-\right| { \left| {c_n^--c \over \mu}\right|^{\rho-1}}  \\
\nonumber &\le&  \eta^{-1} \left|p_n-c_n^-\right| { \left| {c_n^--c \over \mu}\right|^{\rho-1}}  \\
\nonumber &\le& \eta^{-2} \left|p_{n-1}-c_{n-1}^-\right|   \left| {c_n^--c \over \mu} \right|^{\rho-1} { \left| {c_{n-1}^--c \over \mu}\right|^{\rho-1} } \\
\label{mark2} &\le& \eta^{-n} \left|p_1-c_1^- \right|  \prod_{k=1}^n \left| {c_k^--c \over \mu}\right|^{\rho-1} \le \eta^{-n}  e^\pi c_1^- {|L|^\rho \over c^\rho}  \prod_{k=1}^{n} \left| {c_k^--c \over \mu}\right|^{\rho-1}.
\end{eqnarray}

We will now obtain an  estimate on ${(c_k^--c ) / \mu}$.  To that end, first notice, that
$$f_1'(x) \ge { e^{-\pi} (1-c_1^+) \rho \over \mu } \left( {x-c \over \mu }  \right)^{\rho-1} \ge   { e^{-\pi} (1-c_1^+) \rho \over \mu } \left( {\Theta+|R| \over \mu }  \right)^{\rho-1},$$
for all $x \ge f_1^{-1}(c)$, therefore, using the lower bound on $1-c_1^-$ from Lemma $\ref{lem2}$,
\begin{eqnarray}
\nonumber c_k^--c & \le & 1 -\left( \min_{x \ge f_1^{-1}(c)} \left\{f_1'(x) \right\} \right)^{k-1}(1-c_1^-) -c \\
\nonumber & \le & 1-c-\left({ e^{-\pi} (1-c_1^+) \rho \over \mu } \left( {\Theta +|R|\over \mu }  \right)^{\rho-1} \right)^{k-1} (1-c_1^-) \\
\nonumber &\le& \mu - \left( {e^{-2 \pi} \over \eta }\left( {\Theta +|R| \over \mu }  \right)^{\rho-1} \right)^{k-1} {\eta^n \Theta  \over 1-\eta^n} \implies \\
  \nonumber {c_k^--c \over \mu}  &\le& 1 - \left( {e^{-2 \pi} \over \eta }\left( {\Theta + |R| \over \mu }  \right)^{\rho-1} \right)^{k-1} {\eta^n  \Theta \over \mu(1-\eta^n)}  \implies \\
\nonumber \prod_{k=1}^n \left|{c_k^--c \over \mu} \right|^{\rho-1} &\le& \prod_{k=1}^n \left( 1 - \left( {e^{-2 \pi} \over \eta }\left( {\Theta + |R| \over \mu }  \right)^{\rho-1} \right)^{k-1} {\eta^n \Theta \over \mu (1-\eta^n)} \right)^{\rho-1} \\
\nonumber &\le&  \exp\left( \ln \left(\prod_{k=1}^n \left( 1 - \left( {e^{-2 \pi} \over \eta }\left( {\Theta + |R| \over \mu }  \right)^{\rho-1} \right)^{k-1} {\eta^n \Theta \over \mu (1-\eta^n)} \right)^{\rho-1} \right) \right) \\
\nonumber &\le& \exp \left( -K (\rho-1) {\eta^n \Theta \over \mu (1-\eta^n)} \sum_{k=1}^n \left({e^{-2 \pi} \over \eta} {(\Theta+|R|)^{\rho-1} \over \mu^{\rho-1}}  \right)^{k-1}  \right),
\end{eqnarray}
where $K$ is some immaterial constant of order $1$. Finally, $(\ref{mark2})$ becomes
$$L \le \eta^{-n}  e^\pi c_1^- {| L|^\rho \over c^\rho}   \exp \left( -K (\rho-1) {\eta^n \Theta \over \mu(1-\eta^n)} \sum_{k=1}^n \left({e^{-2 \pi} \over \eta} {(\Theta+|R|)^{\rho-1} \over \mu^{\rho-1}}  \right)^{k-1}  \right).$$
which results in the required lower bound for $L$. 

The lower bound for $R$ is obtained in a similar way.
\end{proof}
\medbreak

\bigskip

\section{{\it A priori} bounds}

In this Section we will demonstrate the existence of a relatively compact (in $\cL^0$) set of Lorenz maps which is invariant under the renormalization. This set, $\cK^\pi_\eps$ is defined as in $(\ref{set2})$, i.e. by imposing bounds on the distortion of the diffeomorphic factors of the Lorenz maps, and by restricting the location of the critical point. The invariance of this set will be demonstrated in two steps. 

First we prove renormalization invariance of the distortion bounds in Proposition \ref{inv_distortion}.  Our proof uses the fact that the diffeomorphic factors of the renormalization extend as diffeomorphisms to much larger intervals than those on which their distortion has to be estimated. This allows us to estimate the distortion for the renormalized maps via the Koebe principle. 

At the second step, we prove invariance of the bounds on the critical point in Proposition \ref{inv-crpoint}. This is done with the help of the crucial bounds on the length of the renormalization intervals $L$ and $R$ from Lemma $\ref{lem3}$.

Both Propositions are quite technical and require a simultaneous adjustment of several quantities. Specifically, we are forced to choose $\rho$ not too small in both Proposition, and  $M>N>(\ln 2 /\ln \rho +1+k)  \rho$ in the first Proposition; in the proof of the second Proposition $N<M<2 \rho /(1 + 3 k \ln \rho / \ln 2)$, in both bounds $k \le 1/2$ is a free parameter.

In the last Proposition \ref{apriori} we show that the range of of the length of combinatorics allowed  in Propositions \ref{inv_distortion} and \ref{inv-crpoint} is non-empty. 

In conclusion, we obtain a priori bounds for $\rho$ which are not too small, and for the length of combinatorics which, essentially, lie in the interval $((\ln2 / \ln \rho +1) \rho, 2 \rho)$.

Recall that by Lemma $\ref{renorm_parts}$, the diffeomorphic coefficients of the renormalized map  are
$$\tilde{\phi}=\xi_{C}^{-1} \circ f_1^n \circ \phi  \circ \xi_{\phi^{-1} \circ f_1^{-n}(C)}, \quad  \tilde{\psi}=\xi_{C}^{-1} \circ f_0^m \circ \psi  \circ \xi_{\psi^{-1} \circ f_0^{-m}(C)}.$$

Also, recall the definition of the subset $\cM_{N,M}$ from the Introduction: this is the subset of $\cM$ (monotone types)  of all $\omega$'s such that the length of words in $\omega$  satisfies $N \le |\omega^-| \le  M$ and  $N \le |\omega^+| \le M$.

The next Proposition establishes the conditions for the invariance of the distortion of the coefficients under renormalization.

\begin{prop}\label{inv_distortion} (Invariance of distortion). 
There exist $\hat{\rho}>2$, such that for every $\rho> \hat{\rho}$, $\beta > \ln 2 /\ln \rho$, $0 < \pi < k \ln(\rho)$, $k\le 1/2$, and $M, N \in \field{N}$, $M>N>(\beta +1 +k)  \rho$, if $f \in \cK^\pi_{\rho^{-\beta}} \cap \cL_{\cM_{N,M}}$ then the diffeomorphic factors of the renormalization satisfy 
$${\rm dist} [\tilde{\phi}] < \pi, \  {\rm dist} [\tilde{\psi}] < \pi.$$
\end{prop}
\begin{proof}
We consider the exponential of the distortion of $\tilde{\phi}$ on  $[0,1]$. For any $x,y \in [0,1]$,
\begin{equation} \label{b0}
 {\tilde{\phi}'(x) \over \tilde{\phi}'(y) } = {(f^n_1 \circ  \phi)' ( \xi_{\phi^{-1} \circ f_1^{-n}(C)}(x))) \over (f^n_1 \circ  \phi)' ( \xi_{\phi^{-1} \circ f_1^{-n}(C)}(y)))}. 
\end{equation}
 Recall, that $C=[p,q]$,  and that, by Lemma  $\ref{psiphi_ext_lemma}$, $f^{-n}_1$ and $f_0^{-m}$ are well-defined diffeomorphisms at least on $(c_1^+,1)$  and $(0,c_1^-)$, respectively.  By the Koebe Principle 
$${\left(f^n_1 \circ \phi\right)'( z  ) \over \left( f^n_1 \circ \phi \right)' (w)} \le \left(1+\tau \over \tau \right)^2,$$
where $z,w \in \phi^{-1}(f_1^{-n}(C))$, and
$$\tau=\max\{\tau_1,\tau_2\}, \quad \tau_1={1-q \over q-p}, \quad \tau_2={p-c_1^+ \over q-p}.$$
Similarly,  for $z,w \in \psi^{-1}(f_0^{-m}(C))$,
$${\left( f^m_0 \circ \psi \right)' ( z  ) \over \left( f^m_0 \circ \psi \right)' (w)} \le \left(1+\zeta \over \zeta \right)^2,$$
where 
$$\zeta=\max\{\zeta_1,\zeta_2\}, \quad \zeta_1={p \over q-p}, \quad \zeta_2={c_1^--q \over q-p}.$$

Notice, that $(x-p)/(x-q)$ is a decreasing function of $x$, therefore, $(1+\zeta_2)/\zeta_2=(c_1^--p)/(c_1^--q) \ge  (1+\tau_1)/\tau_1=(1-p)/(1-q)$. Similarly,  $(1+\tau_2)/\tau_2 \ge  (1+\zeta_1)/\zeta_1$. Therefore,

\begin{eqnarray}
\nonumber  \max\{{\rm dist}[\tilde{\phi}],{\rm dist}[\tilde{\psi}] \} &\le&  \max\left\{\left(1+\tau_2 \over \tau_2 \right)^2,\left(1+\zeta_2 \over \zeta_2 \right)^2\right\} \\
\label{b1}  &=&\max\left\{\left({q-c_1^+ \over p-c_1^+} \right)^2,\left({c_1^--p \over c_1^--q} \right)^2\right\} .
\end{eqnarray}

Below we will demonstrate that $(\ref{b1})$ is less than $e^{\pi}$ for sufficiently large $n$ and $m$.

Recall, that $\Delta$ from Lemma $\ref{delt}$ serves as a  lower bound on $p-f_0^{-1}(c)$, while $\Theta$ is a lower bound on $q-f_1^{-1}(c)$. Then, using that  $p-f_0^{-1}(c) < p-c_1^+$, together with the upper bounds on $L$ and $R$ from Lemma $\ref{lem3}$, we get
\begin{eqnarray}
\nonumber{q-c_1^+ \over p-c_1^+} & \le & 1+{q-p \over p-c_1^+} \le 1+{|C| \over \Delta} = 1+{|L| + |R|  \over \Delta} \\
\nonumber &\le& 1+ \left(  \left( { e^{2 \pi} \mu c^{\rho-1}  \over \rho} \right)^{1 \over \rho+1} \left(\mu^\rho e^{2 \pi}  \over |R|^{\rho+1}  \left(\gamma^{-m}-1 \over \gamma^{-1} -1 \right) + \mu^\rho e^\pi |L|   \right)^{\rho-1 \over \rho (\rho+1) }  \left( {1-\gamma  \over  \gamma^{2-n}-\gamma } \right)^{1 \over \rho+1}   + \right.\\
\nonumber &\phantom{\le}&\phantom{1} +  \left.  \left( { e^{2 \pi} c \mu^{\rho-1}  \over \rho} \right)^{1 \over \rho+1} \left(c^\rho e^{2 \pi}  \over |L|^{\rho+1}  \left(\gamma^{-n}-1 \over \gamma^{-1} -1 \right) + c^\rho e^\pi |R|   \right)^{\rho-1 \over \rho (\rho+1) }  \left( {1-\gamma  \over  \gamma^{2-m}-\gamma } \right)^{1 \over \rho+1}  \right)   \times \\
\nonumber &\phantom{\le}& \phantom{1} \times  \left({   c_1^-  \rho \over c}  {\left( 1-{c_1^+ \over c}  \right)^{\rho-1}   \over  \nu^{\rho \over \rho-1} e^{-\pi {\rho \over \rho -1} }  }  \right)^{\rho^n \over \rho^n-1} \\
\nonumber &\le& 1+ \left(  \left( { e^{2 \pi} \mu c^{\rho-1}  \over \rho} \right)^{1 \over \rho+1} \left( e^{\pi}  \over |L|   \right)^{\rho-1 \over \rho (\rho+1) }  \left( {1-\gamma  \over  \gamma^{2-n}-\gamma } \right)^{1 \over \rho+1}   + \right.\\
\nonumber &\phantom{\le}&\phantom{1} +  \left.  \left( { e^{2 \pi} c \mu^{\rho-1}  \over \rho} \right)^{1 \over \rho+1} \left(e^{\pi}  \over |R|   \right)^{\rho-1 \over \rho (\rho+1) }  \left( {1-\gamma  \over  \gamma^{2-m}-\gamma } \right)^{1 \over \rho+1}  \right)   \left({   c_1^-  \rho \over c}  {\left( 1-{c_1^+ \over c}  \right)^{\rho-1}   \over  \nu^{\rho \over \rho-1} e^{-\pi {\rho \over \rho -1} }  }  \right)^{\rho^n \over \rho^n-1}.
\end{eqnarray}
Now, let $N=\min\{n,m\}$, and $M=KN=\max\{n,m\}$, and suppose that $N$ is large, then 
\begin{equation}
 \label{sufcond1} {q-c_1^+ \over p-c_1^+}  \le 1+ 2   e^{\pi p_1} (1-\vareps)^{p_2} \vareps^{p_3}\rho^{p_4},
\end{equation}
where 
\begin{eqnarray}
\nonumber p_1&=& { {K N +1 \over \rho (\rho+1)} + {\rho-1 \over \rho (\rho+1)  } +{2 \over \rho+1} +2 {K N-2 \over \rho+1} +{\rho^{N+1} \over (\rho-1) (\rho^N-1)} },\\
\nonumber p_2&=&{\rho-1 \over \rho+1} -{1 \over \rho+1} -{\rho^N \over \rho^N-1}, \\
\nonumber p_3&=&-{K N \over \rho (\rho+1)}+{1 \over \rho+1} -{\rho^{N+1} \over (\rho-1) (\rho^N-1)}, \\
 \nonumber p_4&=& -{1 \over \rho+1}-{N -2 \over \rho+1 }+{K N \over \rho (\rho+1)} +{\rho^N \over \rho^N-1}.
\end{eqnarray} 
In particular, a sufficient condition for $(q-c_1^+) / (p-c_1^+)$ to be  smaller than $e^{\pi/2}$ is 
\begin{equation}
\label{sufcond2}1+ 2   e^{\pi p_1} (1-\vareps)^{p_2} \vareps^{p_3}\rho^{p_4} < e^{\pi \over 2}.
\end{equation}
Set $\vareps=\rho^{-\beta}$ for some $\ln 2 / \ln \rho < \beta$ (this ensures that $\vareps < 1/2$). Then equation $(\ref{sufcond2})$ becomes 
\begin{equation}
\label{sufcond3}1+ {\rm const} \  e^{\pi p_1} \rho^{q_2} < e^{\pi \over 2}.
\end{equation}
where the constant, which is an upper bound on $2 (1-\rho^{-\beta})^{p_2}$, $\rho>2$, can be chosen to be independent of $N$, $\rho$ and $K$. 
\begin{eqnarray}
\nonumber q_2&=&p_4-\beta p_3 \\
\nonumber &=&  { ( ( (\beta + 1) \rho - \beta - 1)  K - (\rho^2 - \rho)) N + ((\beta + 1) \rho^3 + (\beta - 2) \rho + \rho^2)\over (\rho^3 - \rho) \rho^N -\rho^2 + \rho} \rho^N -\\
\nonumber &&-{ (\beta - 1) \rho^2 + ((- (\beta + 1) \rho + \beta + 1) K + \rho^2 - \rho) N - (\beta - 1) \rho \over \rho^3 - (\rho^3 - \rho) \rho^N - \rho}.
\end{eqnarray}
Let $N=t \rho$ and $K=1$. Consider the product of the two factors $e^{\pi p_1}  \rho^{q_2}$ with the upper bound on the distortion $k \ln \rho$ in place of $\pi$: $\rho^{k p_1} \rho^{q_2}$. The power of $\rho$ in this expression has the following limit:
$$\lim_{\rho \rightarrow \infty} k p_1 +q_2=-t+\beta+k+1.$$
For this choice of $\pi$, $K$ and $N$,  the condition $(\ref{sufcond3})$ for large $\rho$ takes the form
$$1+ {\rm const} \rho^{-t+\beta+k+1} <\rho^{k/2},$$ 
and  can be satisfied if $t>\beta+k+1$.

In a similar way $(c_1^--p) / (c_1^+-q)$ is less than $e^{\pi /2}$ for the choice of $\pi$, $N$, $M$ and $\rho$ as in the hypothesis.
\end{proof}
\medbreak

Recall that according to Lemma $\ref{renorm_parts}$, the critical point of a renormalized Lorenz map is given by 
\begin{equation}\label{tildec}
\tilde{c}={|L| \over |C|}.
\end{equation}

\begin{prop}(Invariance of the bounds on the critical point). \label{inv-crpoint}
There exist $\hat{\rho}>2$, such that for every $\rho> \hat{\rho}$, $\beta > \ln 2 /\ln \rho$, $0 < \pi < k \ln(\rho)$, $k\le 1/2$, and $M, N \in \field{N}$, $N<M<2 \beta \rho /(\beta + 3 k)$, if $f \in \cK^\pi_{\rho^{-\beta}} \cap \cL_{\cM_{N,M}}$ then the critical point $\tilde{c}$ of the renormalization satisfies 
$$\tilde{c} \in (\rho^{-\beta},1-\rho^{-\beta}).$$   
\end{prop}
\begin{proof} We will start with the lower bound on $\tilde{c}$. 

According to $(\ref{tildec})$, for $\tilde{c}$ to be larger or equal to some $\vareps>0$ it is sufficient that
\begin{eqnarray}
\label{epsinv} {1 \over 1+\max {|R| \over |L|} } \ge \vareps \Leftrightarrow 1 \ge \vareps \left( 1+ \max {|R| \over |L|} \right).
\end{eqnarray}
The maximum of the ratio of the lengths of $R$ and $L$ can be estimated using bounds from Lemma $\ref{lem3}$:
\begin{eqnarray}
\nonumber {|R| \over |L|} &\le& \left( {c e^{2 \pi} \mu^{\rho-1}  \over \rho} \right)^{1 \over \rho+1} \left(c^\rho e^{2 \pi} \left(\gamma^{-1} -1 \over \gamma^{-n}-1  \right)   \right)^{\rho-1 \over \rho (\rho+1)} {1 \over |L|^{\rho-1 \over \rho}} \left( { 1-\gamma \over \gamma^{2-m}-\gamma  } \right)^{1 \over \rho+1}     \left( 1 \over 1+{e^{-\pi} |R| \over |L|^{\rho+1} \left({ \gamma^{-n}-1 \over \gamma^{-1}-1}   \right) }  \right)^{\rho-1 \over \rho (\rho+1)}     \over |L| \\
\nonumber & =&  {\left( {c e^{2 \pi} \mu^{\rho-1}  \over \rho} \right)^{1 \over \rho+1} \left(c^\rho e^{2 \pi} \left(\gamma^{-1} -1 \over \gamma^{-n}-1  \right)   \right)^{\rho-1 \over \rho (\rho+1)} \left( { 1-\gamma \over \gamma^{2-m}-\gamma  } \right)^{1 \over \rho+1}  \left( 1 \over 1+{e^{-\pi} |R| \over |L|^{\rho+1} \left({ \gamma^{-n}-1 \over \gamma^{-1}-1}   \right) }  \right)^{\rho-1 \over \rho (\rho+1)} \over |L|^{2 \rho -1 \over \rho} }.
\end{eqnarray}
We use the bounds from Lemma $\ref{lem3}$ again:
\begin{eqnarray}
\nonumber &&{|R| \over |L|}  \left( 1+{e^{-\pi} |R| \over |L|^{\rho+1} \left({ \gamma^{-n}-1 \over \gamma^{-1}-1}  \right) }   \right)^{\rho-1 \over \rho (\rho+1)}   \le  { \left( {c e^{2 \pi} \mu^{\rho-1}  \over \rho} \right)^{1 \over \rho+1} \left(c^\rho e^{2 \pi} \left(\gamma^{-1} -1 \over \gamma^{-n}-1  \right)   \right)^{\rho-1 \over \rho (\rho+1)} \left( { 1-\gamma \over \gamma^{2-m}-\gamma  } \right)^{1 \over \rho+1} \over \left( {e^{-\pi} c^\rho \over  c_1^-}  \left( {e^{-\pi} \mu \over (1-c_1^+) \rho} \right)^n \right)^{{2 \rho -1 \over \rho(\rho-1)}} }  \implies \\
\nonumber &&\left( {|R| \over |L|} \right)^{\rho (\rho+1) \over \rho-1}  \left( 1+{e^{-\pi} |R| \over |L|^{\rho+1} \left({ \gamma^{-n}-1 \over \gamma^{-1}-1}  \right) }   \right)  \le  { \left( {c e^{2 \pi} \mu^{\rho-1}  \over \rho} \right)^{\rho \over \rho-1} \left(c^\rho e^{2 \pi} \left(\gamma^{-1} -1 \over \gamma^{-n}-1  \right)   \right) \left( { 1-\gamma \over \gamma^{2-m}-\gamma  } \right)^{\rho \over \rho-1} \over \left( {e^{-\pi} c^\rho \over  c_1^-}  \left( {e^{-\pi} \mu \over (1-c_1^+) \rho} \right)^n \right)^{{ (2 \rho -1) (\rho+1) \over (\rho-1)^2}} }  \implies \\
\nonumber &&\left({|R| \over |L|}\right)^{2 \rho^2+ \rho-1 \over \rho-1   } {e^{-\pi} \over |R|^{\rho} \left({ \gamma^{-n}-1 \over \gamma^{-1}-1}  \right) }    \le  { \left( {c e^{2 \pi} \mu^{\rho-1}  \over \rho} \right)^{\rho \over \rho-1} \left(c^\rho e^{2 \pi} \left(\gamma^{-1} -1 \over \gamma^{-n}-1  \right)   \right) \left( { 1-\gamma \over \gamma^{2-m}-\gamma  } \right)^{\rho \over \rho-1} \over \left( {e^{-\pi} c^\rho \over  c_1^-}  \left( {e^{-\pi} \mu \over (1-c_1^+) \rho} \right)^n \right)^{{ (2 \rho -1) (\rho+1) \over (\rho-1)^2}} }  \implies  \\
\nonumber &&\left({|R| \over |L|}\right)^{2 \rho^2+ \rho-1 \over \rho-1   } \hspace{-9mm} {e^{-\pi} \over \left(c \mu^{\rho-1} e^\pi   \right)^{\rho \over \rho+1} \gamma^{(m-1){\rho \over \rho+1} -(n-1)}}   \le  { \left( {c e^{2 \pi} \mu^{\rho-1}  \over \rho} \right)^{\rho \over \rho-1} \left(c^\rho e^{2 \pi} \left(\gamma^{-1} -1 \over \gamma^{-n}-1  \right)   \right) \left( { 1-\gamma \over \gamma^{2-m}-\gamma  } \right)^{\rho \over \rho-1} \over \left( {e^{-\pi} c^\rho \over  c_1^-}  \left( {e^{-\pi} \mu \over (1-c_1^+) \rho} \right)^n \right)^{{ (2 \rho -1) (\rho+1) \over (\rho-1)^2}} }  \implies  \\
\label{eqRL2}  &&{|R| \over |L|}  \le  e^{\pi q_1} c^{q_2}  \mu^{q_3} \rho^{q_4},
\end{eqnarray}
where
\begin{eqnarray}
\nonumber q_1&=& { 2 \rho^3 + (2 \rho^3 + (4 \rho^3 + 3 \rho^2 - 2 \rho - 1) K - 4 \rho^2 + 2 \rho) N + 2 \rho^2 - 2 \rho + 2 \over 2 \rho^4 + \rho^3 - 3 \rho^2 - \rho + 1 },\\
\nonumber q_2&=& -{ \rho^4 + 2 \rho^3 + 3 \rho^2 - 2 \rho \over 2 \rho^4 + \rho^3 - 3 \rho^2 - \rho + 1 },\\
\nonumber q_3&=&{2 \rho^4 - (2 \rho^3 + 3 \rho^2 - 1) K N - 4 \rho^3 + 2 \rho^2 \over 2 \rho^4 + \rho^3 - 3 \rho^2 - \rho + 1},\\
\nonumber q_4&=& { 2 \rho^3 - (\rho^3 - (\rho^3 + 5 \rho^2 - \rho - 1) K - \rho) N - 2 \rho^2 \over 2 \rho^4 + \rho^3 - 3 \rho^2 - \rho + 1 }
\end{eqnarray}
Next, set $\vareps=\rho^{-\beta}$ for some $\beta>\ln(2)/\ln(\rho)$ (this guarantees that $\vareps < 1/2$). Consider the product of the three factors $e^{\pi q_1}  \rho^{-\beta (q_3+1)} \rho^{q_4}$ with the upper bound on the distortion $k \ln \rho$ in place of $\pi$: $\rho^{k q_1} \rho^{-\beta (q_3+1)} \rho^{q_4}$. The power of $\rho$ in this expression is as follows
\begin{eqnarray}
\nonumber k q_1+q_4-\beta (q_3+1)\!\!&\!\!=\!\!&\!\!{ 2  (1+ 2 K) \rho^3  N \over 2 \rho^4 + \rho^3 - 3 \rho^2 - \rho + 1 }-{4 b \rho^4 +(1 - K (2 b + 1)) \rho^3 N \over  2 \rho^4 + \rho^3 - 3 \rho^2 - \rho + 1 }+ \\
\nonumber &\!\!\phantom{\!\!=\!\!}\!\!&  + {(3 b + 2) \rho^3 +(b - 2) \rho^2 +(((3 b + 5) \rho^2 - b - \rho - 1) K + \rho) N + b \rho - b \over 2 \rho^4 + \rho^3 - 3 \rho^2 - \rho + 1 } \\
\label{power} &\!\!\phantom{\!\!=\!\!}\!\!&+{ 2 \rho^3 + ((3 \rho^2 - 2 \rho - 1) K - 4 \rho^2 + 2 \rho) N + 2 \rho^2 - 2 \rho + 2 \over 2 \rho^4 + \rho^3 - 3 \rho^2 - \rho + 1 },
\end{eqnarray}
and has the following limit in the case $N= t \rho$, $K=1$ (the leading terms are written in the first line of $(\ref{power})$):
$$\lim_{\rho \rightarrow \infty}  {q_1 \over 2} +q_4-\beta( q_3+1)=(t-2) \beta +3 k t.$$

Therefore, for sufficiently large $\rho$, there exists $K>1$ and $N< 2 \beta \rho/(\beta + 3 k)$ such that $(\ref{power})$ is negative, and the inequality 
$$e^{\pi q_1} (1-\rho^{-\beta})^{q_2}  \rho^{-\beta q_3} \rho^{q_4} < (1-\rho^{-\beta}) \rho^\beta$$
is satisfied for all $0<\pi < k \ln \rho$.

Similarly, $ {|L| \over |R|} < (1-\vareps)/\vareps$ for $\rho$ and $\pi$ as in the hypothesis.
\end{proof}
\medbreak

The following results is an almost immediate corollary of Propositions $\ref{inv_distortion}$ and $\ref{inv-crpoint}$. The only thing remaining to be verified, is that the bounds on $M$ and  $N$ from Propositions  $\ref{inv_distortion}$ and $\ref{inv-crpoint}$ indeed specify a non-empty range of combinatorics. For that, it is sufficient that 
$$(\beta +k+1) \rho < {2 \beta \rho \over \beta + 3 k} \implies k \in \left(0 ,  {\sqrt{(4 \beta +3)-12 (\beta^2-\beta)} -4 \beta -3 \over 6} \right),$$
the last interval being indeed non-empty.

\begin{prop}(A priori bounds). \label{apriori}
There exist $\hat{\rho}>2$, such that for every $\rho> \hat{\rho}$, $\ln 2 / \ln \rho <\beta<1$, $0 < \pi < k \ln \rho$, 
$$k < {\sqrt{4 \beta^2+36 \beta +9} -4 \beta -3 \over 6},$$ 
and $M, N \in \field{N}$, 
$$(\beta +k+1) \rho < N \le  M < {2 \beta \rho \over \beta + 3 k},$$ 
one has
$$\cR[\cK^{\pi}_{\rho^{-\beta}} \cap \cL_{\cM_{N,M}}^S] \subset \cK^\pi_{\rho^{-\beta}}.$$
\end{prop}

\bigskip

\section{Periodic points of renormalization}

We consider a restriction $\cR_\omega$ of the renormalization operator to some 
$$\omega=(0 \overbrace{1 \ldots 1}^{n}, 1 \overbrace{0 \ldots 0}^{m})  \in \cM, \quad N < n < M, \ N < m  < M,$$
where $N$ and $M$ are as in Proposition $\ref{apriori}$.

In this Section we will demonstrate that $\cR_\omega$ has a fixed point. We will generally follow the approach  of  \cite{Win2} (and we will make a conscientious attempt to keep the notation in line with that work). One important difference with the case considered in \cite{Win2}, however, is that we are looking at a  different class of return times. This will introduce some extra difficulties, especially evident in the proof of Lemma $\ref{iff_lemma}$, somewhat more involved than its analogue Lemma $6.8$  from \cite{Win2}.

We will start by quoting several previously established results.

\begin{definition}
A branch $I$ of $f^n$ is {\it full} if $f^n$ maps $I$ onto the domain of $f$. $I$ is {\it trivial} if $f^n$ fixes both  endpoints of $I$.
\end{definition}

We  will now quote several facts about Lorenz  maps, established in \cite{MM}.

\begin{definition}
A  slice in the parameter plane is any set of the form 
$$\cS=[0,1]^2 \times \{c\} \times \{\phi\} \times \{\psi\},$$
where $c$, $\phi$  and $\psi$ are fixed. We will use the simplified notation $(u,v) \in \cS$.
\end{definition}

A slice $\cS$ induces a family of Lorenz maps  
$$\cS \ni (u,v) \mapsto  (u,v,c,\phi,\psi) \subset \cL^0.$$
any family induced by a  slice  is {\it full}, that is it contains maps of all possible combinatorics. Specifically (see \cite{MM}  for details), 

\begin{prop}\label{theoremA}({\it Theorem A from \cite{MM}}).
Let $(u,v) \mapsto (u,v,c,\phi,\psi)$ be a family induced by a slice. Then this family intersects $\cL_{\bar{\omega}}^0$ for every $\bar{\omega}$ (finite or infinite) such  that $\cL_{\bar{\omega}}^0  \ne \emptyset$.  
\end{prop}

\begin{lemma}\label{lemma41}({\it Lemma 4.1 from  \cite{MM}}).
Assume that $f$ is renormalizable. Let $(l,c) \supset L$ be the branch of $f^{n+1}$ and $(c,r) \supset  R$  be that  of $f^{m+1}$. Then
$$f^{n+1}(l) \le l, \quad f^{m+1}(r) \ge r.$$
\end{lemma}

Let $\pi$, $\vareps$ and  $\cK^{\pi}_{\vareps}$ be as in the previous Section.  Consider the set 
\begin{equation}
\cY = \cL_\omega^S \cap \cK^{\pi}_{\vareps}. 
\end{equation}

\begin{prop}
The boundary of $\cY$ consists of three parts: $f \in \partial \cY$ iff at least one of the following holds:
\begin{itemize}
\item[C1.] the left and the right branches of $\cR[f]$ are full or trivial;
\item[C2.] ${\rm dist}[\phi]=\pi$ or ${\rm dist}[\psi]=\pi$;
 \item[C3.] $c(f)=\vareps$, or $c(f)=1-\vareps$.
\end{itemize}
\end{prop}
\begin{proof}
Consider the boundary of $\cL_\omega^0$. If either branch of $\cR_\omega[f]$ is full or trivial, then there exists an perturbation of $f$,  however small, such that $f$ is no longer renormalizable. Hence C1 holds on $\partial \cL_\omega^0$. If $f \in \cL_\omega^0$ does not satisfy C1 then, according to Lemma $\ref{lemma41}$,  all small perturbations of it will be still renormalizable.

Conditions C2 and C3 are part of the boundary of $\cK^{\pi}_{\vareps}$. By Proposition $\ref{theoremA}$ these boundaries intersect $\cL_\omega^S$, and hence C2 and C3 are also the boundary conditions for $\cY$.  
\end{proof}
\medbreak

Fix $c_0 \in (\vareps,1-\vareps)$, and let $\cS=[0,1]^2 \times \{c_0\} \times \{\rm id \} \times \{\rm id\}$.  Recall, that the linear structure on the space $\cD^2$ is defined via the nonlinearity operator:
$$ \alpha \phi+\beta \psi=  N^{-1}_{\left(\alpha N_{\phi}+\beta N_{\psi}  \right)}.$$

Introduce the {\it deformation retract} onto $\cS$ as 
\begin{eqnarray}
\nonumber \pi_t (u,v,c,\phi,\psi)&\equiv &(u,v,c+t(c_0-c),\phi_t, \psi_t) \\
\label{retract} &=&  (u,v,c+t(c_0-c), (1-t) \phi+t \ {\rm id}, (1-t)  \psi+t \ {\rm id}).
\end{eqnarray}

Let 
$$\cR_t=\pi_t \circ \cR.$$

We will strengthen the conditions on the set $\cY$ and consider a smaller set
\begin{equation}\label{cYd}
\cY_\delta = \cY \cap \{f \in \cL^S_\omega: c(\cR[f])  \ge  \delta\}. 
\end{equation}

The boundary of $\cY_\delta$ is given by conditions $C1$-$C3$ together with 
\begin{itemize}
\item[C4.] $\{f \in \cY: c(\cR[f])=\delta\}$.
\end{itemize}

\begin{lemma}\label{iff_lemma}
There exists a choice of $c_0$ in $(\ref{retract})$ and $\delta \in (\vareps,1-\vareps) $, such that $\cR$  has a fixed point in $\partial  \cY_\delta$ iff $\cR_t$ has a fixed point in $\partial \cY_\delta$ for some $t \in [0,1]$.
\end{lemma}
\begin{proof}
The direct statement is obvious since $\cR \equiv \cR_0$. 

The converse is also obvious when $t=0$, and we, therefore, consider $t>0$. Assume that $f \in \partial \cY_\delta$ with the coefficients $(\phi,\psi)$ is such that $\cR_t f=f$ for some $t \in (0,1]$, and  assume that $\cR$ has no fixed point on  $\partial \cY_\delta$. We  will demonstrate that this is impossible.

Choose $c_0$ close to $1-\vareps$: $c_0=1- \vareps -\nu$ for some small $\nu$. By Proposition $\ref{inv-crpoint}$, $c({\cR[f]}) \in [\vareps,1-\vareps]$ whenever $c(f)$ is. Together with the condition $c(\cR[f])  \ge  \delta$ this implies that
$$c({\cR[f]}) \in [\delta,1-\vareps].$$

Since $t>0$, by formula $(\ref{retract})$ $c(\cR_t[f])$ is strictly in the interior of $[\delta,1-\vareps]$ for all $t \in (0,1]$.  Therefore, neither C3 nor C4 can hold for $f=\cR_t[f]$ for $t \in (0,1]$.

The distortion of the coefficients of $\cR[f]$ is not  greater  than $\pi$ by Proposition $\ref{inv_distortion}$. For $t \in (0,1]$  distortion of the diffeomorphic parts  $(\tilde{\phi}_t,\tilde{\psi}_t)$  of $\cR_t[f]$ is strictly smaller than that of $(\tilde{\phi},\tilde{\psi})$ (diffeomorphic coefficients for $\cR[f]$). This can be seen from the following computation:
$${\tilde{\phi}_t'(x) \over \tilde{\phi}_t'(y) }={ {\exp{[\int_0^y (1-t)  N_{\tilde{\phi}}(s) + t N_{{\rm id}}(s) d s]} \over \int_0^1 {  \exp{[ \int_0^r (1-t)  N_{\tilde{\phi}}(s) + t N_{{\rm id}}(s) d s]}  d r } } \over  {\exp{[\int_0^x (1-t)  N_{\tilde{\phi}}(s) +t N_{{\rm id}}(s) d s]} \over \int_0^1 {  \exp{[ \int_0^r (1-t)  N_{\tilde{\phi}}(s) +t N_{{\rm id}}(s) d s]}  d r } } }={ {\exp{[\int_0^y (1-t)  N_{\tilde{\phi}}(s)  d s]}} \over  {\exp{[\int_0^x (1-t)  N_{\tilde{\phi}}(s) d s]}} } =\left(\tilde{\phi}'(y) \over \tilde{\phi}'(x)  \right)^{1-t} < e^\pi.
$$
Similarly for $\tilde{\psi}_t$.  Therefore, we have that C2 does not hold for $f=\cR_t[f]$ for $t \in (0,1]$.

The only possibility is that, if $f=\cR_t[f] \in \partial \cY_\delta$ then it belongs to the part of the boundary described by C1.

Suppose that  either branch of $\cR[f]$ is full; for definitiveness, suppose $c_1^-(\cR[f])=1$. Since $\phi$ fixes both end points of the unit interval, this implies that $u(\cR[f])=1$, and since the deformation retract does not change the value of $u$, $u(\cR_t[f])=1$. Since $\phi_t$ fixes $1$ as well, we get that $c_1^-(\cR_t[f])=1$, and therefore, the corresponding branch of $\cR_t[f]$ is full as well.  This shows that $f$ can not be fixed by $\cR_t$ since a renormalizable map  can not have a full branch. Therefore, one of the branches of $\cR[f]$ must be trivial.

Before we proceed with the last case of trivial branches, we will derive an upper bound on  $\phi_t(u)$ and a lower bound on $\psi_t(1-v)$. Recall, that $\phi_t=(1-t) \phi + t \ {\rm id}$ where the linear structure is given by $(\ref{lin_struct})$. Then, on one hand,
\begin{eqnarray}
\nonumber \phi_t(x)&= &{\int_0^x \left(\phi'(r)\right)^{1-t}  dr \over \int_0^1  \left(\phi'(r)\right)^{1-t} d r } = 1 -{\int_x^1 \left(\phi'(r)\right)^{1-t}  dr \over \int_0^1  \left(\phi'(r)\right)^{1-t} d r }  \le  1 -{\int_x^1 \left(\phi'(r)\right)^{1-t}  dr \over \left( \int_0^1  \phi'(r) d r \right)^{1-t} }  \\
 \nonumber  &=& 1 -\int_x^1 \left(\phi'(r)\right)^{1-t}  dr  = \int_0^x \left(\phi'(r)\right)^{1-t}  dr \le   \int_0^x  \left({c_1^- \over u} e^\pi  \right)^{1-t} d r \le  \left({c_1^- \over u} e^\pi  \right)^{1-t} x, 
\end{eqnarray}
and 
\begin{equation} \label{phitu1}
\phi_t(u) \le (c_1^-)^{1-t} u^t e^{\pi (1-t)}.
\end{equation}
On the other hand,
$$\phi_t(x) \le\int_0^x \left(\phi'(r)\right)^{1-t}  dr \le  \int_0^x \phi'(r)  dr \sup_{r \in (0,x)} \left(\phi'(r) \right)^{-t} \le  \phi(x) \left(e^\pi {x \over \phi(x)} \right)^t,$$
and
\begin{equation} \label{phitu2}
\phi_t(u) \le (c_1^-)^{1-t} u^t e^{\pi t}.
\end{equation}
We can now take a linear combination of $(\ref{phitu1})$ and  $(\ref{phitu2})$ as an upper bound on $\phi_t(u)$. A particularly convenient choice is 
\begin{equation} \label{phitu}
\phi_t(u) \le (c_1^-)^{1-t} u^t \left( t e^{\pi (1-t)} +(1-t) e^{\pi t}\right),
\end{equation}
Notice, that the maximum of the function $\left( t e^{\pi (1-t)} +(1-t) e^{\pi t}\right)$ is achieved at $t=1/2$. 

In a similar way,
\begin{equation} \label{psitv}
\psi_t(1-v) \ge (c_1^+)^{1-t} (1-v)^t \left( t e^{-\pi (1-t)} +(1-t) e^{-\pi t}\right).
\end{equation}

Suppose, the left branch is trivial: $c(\cR[f]) \ge c^-_1(\cR[f])$. Recall, that for a map renormalizable with monotone combinatorics,  $c_1^->f_1^{-1}(c)$, and according to Lemmas $\ref{lem3}$ and $\ref{delt}$  the differences $c_1^--c>|R|+\Theta \ge K$ and $c-c_1^+>|L|+\Delta \ge J$,  where  $K$ and $J$  depend on $\rho$, $\pi$, $\vareps$, $n$ and $m$, but do not depend on the particular form of the map. Also, for large $n$ and $m$, $K$ and $J$ become independent of $n$ and $m$:
$$J=O\left( \kappa {\left({c \over c-c_1^+ } \right)^{\rho-1}}  \left( \nu^{\rho \over \rho-1}  e^{-\pi  \over \rho-1} \right)\right), \quad K=O\left(\eta {\left({\mu \over c_1^--c } \right)^{\rho-1}}  \left( \xi^{\rho \over \rho-1}  e^{-\pi  \over \rho-1} \right)  \right),$$
which follows from the expressions $(\ref{pf0})$ for $\Delta$ and $(\ref{qf1})$ for $\Theta$.

Next, suppose $\nu$ is small: $\nu << K$. Then, on one hand, 
$$c_1^-(\cR_t[f])-c(\cR_t[f]) \le (c_1^-(\cR[f]))^{1-t} u(\cR[f])^t \left( t e^{\pi (1-t)} +(1-t) e^{\pi t}\right) -c(\cR[f])- t (c_0-c(\cR[f])).$$
Recall, that by Proposition $\ref{inv_distortion}$, $\pi$ can be chosen small if one considers large $n$ and $m$. Since $c_1^-(\cR[f])-c(\cR[f])<0$, the expression  
$$(c_1^-(\cR[f]))^{1-t} u(\cR[f])^t \left( t e^{\pi (1-t)} +(1-t) e^{\pi t}\right) -c(\cR[f]) <K-\nu$$ 
if $\pi$ is small. Then  since $c_0-c(\cR[f])$ is larger than $-\nu$, we have that
$$(c_1^-(\cR[f]))^{1-t} u(\cR[f])^t \left( t e^{\pi (1-t)} +(1-t) e^{\pi t}\right) -c(\cR[f])- t (c_0-c(\cR[f]))< K,$$
and hence $f=\cR_t[f]$ is not renormalizable with the monotone combinatorics $\omega$.

Now, suppose that the right branch is trivial. Then
\begin{eqnarray}
\nonumber c(\cR_t[f])-c_1^+(\cR_t[f]) \le c(\cR[f]) &+& t(c_0-c(\cR[f])) \\
\nonumber &-&  (c_1^+(\cR[f]))^{1-t} (1-v(\cR[f]))^t \left( t e^{-\pi (1-t)} +(1-t) e^{-\pi t}\right).
\end{eqnarray}

Since $c(\cR[f]) -  c_1^+(\cR[f])<0$, we have that for a sufficiently small $\pi$, 
$$c(\cR[f])-  (c_1^+(\cR[f]))^{1-t} (1-v(\cR[f]))^t \left( t e^{-\pi (1-t)} +(1-t) e^{-\pi t}\right) < {J \over 2},$$
while
$$c(\cR[f]) \!- \! (c_1^+(\cR[f]))^{1-t} (1\!-\! v(\cR[f]))^t \left( t e^{-\pi (1-t)} \!+\!(1\!-\!t) e^{-\pi t}\right)  \!+\! t(c_0-c(\cR[f]))\le {J \over 2} \!+\!(1\!-\!\vareps\!-\!\nu\!-\!\delta).$$
Therefore, the map $\cR_t[f]$ is not renormalizable with the monotone combinatorics $\omega$ for  $t \in [0,1]$, if we chose $\delta$ so that
\begin{equation} \label{Jineq}
1-\vareps-\nu-\delta < J/2. 
\end{equation}

We now notice, that according to Lemmas  $\ref{lem3}$ and $\ref{delt}$
\begin{equation} \label{growthrate}
J=O\left( \left(\delta (1-\delta)^{\rho \over \rho-1} \right)^{\rho^n \over \rho^n-1} \right)+O\left(\delta^{\rho \over \rho-1} (1-\delta)^{n \over \rho-1} \right) \exp\left( O \left( (1-\delta)^{{\rho^m \over \rho^m-1}-{\rho-1 \over \rho+1} (n-1) }  \right)  \right).
\end{equation}

If $\delta$ is small, then the above expression demonstrates that $J=O \left(\delta^{\rho^n \over \rho-1} \right)$, and the inequality $(\ref{Jineq})$ is not satisfied. On the other hand, if $\delta$ is close to $1-\varepsilon-\nu$, then the exponential in $(\ref{growthrate})$  becomes large and dominates others terms, and $(\ref{Jineq})$ is easily satisfied. Therefore, there exists $\delta \in (0,1-\vareps-\nu)$, not necessarily very close to $1-\vareps-\nu$, such that $(\ref{Jineq})$ holds for all $c > \delta$. 

We conclude that $f=\cR_t[f] \notin \partial \cY_\delta$ which is a contradiction with the assumption in the beginning of the proof.
\end{proof}
\medbreak

According to Theorem B in \cite{MM} the intersection of $\cS$ with $\cL_\omega^S$ contains a connected component $I$ of the interior, called a {\it full island}, such that the family $I \ni (u,v) \mapsto \cR[f]$ is full.

\begin{lemma}\label{ext_lemma}
Any extension of $\cR_1 \arrowvert_{\partial \cY_\delta}$ to $\cY_\delta$ has a fixed point. 
\end{lemma}
\begin{proof}
Assume that $\cR_1$ has no fixed point in  $\partial \cY_\delta$ (otherwise the theorem is trivial). 

Let  $\cS=[0,1]^2 \times \{c_0\} \times \{{\rm id}\} \times \{{\rm id}\}$, where $c_0$ is as in the previous Lemma. This set contains a full island $I$ with $\partial I \subset \partial \cY_\delta$.

Pick any $R: I \mapsto \cS$ such that $R \arrowvert_{\partial I}=\cR_1 \arrowvert_{\partial  I}$. Define the displacement map $d: \partial I \mapsto \field{T}^1$ by
$$d(x)={x-R(x)  \over |x-R(x)|},$$
which is well-defined since $R$ does not have fixed points on $\partial I \subset \partial \cY_\delta$. The degree of $d$ is non-zero since $I$  is full. Therefore, $R$ has a fixed point in $I$ (otherwise $d$ would extend to all of $I$, and would have a degree zero). 
\end{proof}
\medbreak

To finish the proof  of the existence of the fixed points we will require the following theorem from \cite{GD}:

\begin{thm}\label{TheoremGD}
Let $X \subset Y$ where $X$ is closed and $Y$ is a normal topological space. If $f: X \mapsto Y$ is homotopic to a map $g: X \mapsto Y$ with the property that every extension of $g\arrowvert_{\partial X}$ to $X$ has a fixed point in $X$, and if the homotopy $h_t$ has no fixed point on $\partial X$ for every $t \in [0,1]$, then $f$ has a fixed point in $X$. 
\end{thm}

\begin{prop}
$\cR_\omega$ has a fixed point.
\end{prop}
\begin{proof}
$\cR_1$ either has a fixed point in $\partial \cY_\delta$, or otherwise by Lemma $\ref{ext_lemma}$ any of  extensions of $\cR_1 \arrowvert_{\partial \cY_\delta}$ to $\cY_\delta$ has a fixed point. In the second case we can apply  Theorem  $\ref{TheoremGD}$ and Lemma $\ref{iff_lemma}$, to immediately obtain the required result.
\end{proof}
\medbreak

Now we can finish the proof of Theorem A. 

\bigskip

{\it Proof of Theorem A.}
Suppose that $\rho$ is sufficiently large and  $N$ and $M$ are as in Propositions $\ref{inv_distortion}$ and $\ref{inv-crpoint}$. Pick a sequence $\bar{\omega}=(\omega_0, \omega_1, \ldots, \omega_{k-1}), \quad \omega_j \in \cM_{N,M}$. One can use $\cR_{\omega_{k-1}} \circ \ldots \circ \cR_{\omega_0}$ in place of $\cR_\omega$ in the previous Proposition to demonstrate that $\cR_{\omega_{k-1}} \circ \ldots \circ \cR_{\omega_0}$ has a fixed point, which, hence, is a periodic point of $\cR$ of combinatorial type $\bar{\omega}$.
\qed
\medbreak


\bigskip

\section{Appendix}
To present a proof of Theorem B we will need to introduce the concept of transfer maps.

\begin{definition}
An interval $C$ is called a nice interval of $f$, if $C$ is open, the critical point of $f$ is in $C$, and the orbit of the boundary of $C$ is disjoint from $C$. 
\end{definition}

\begin{definition}
Fix $f$ and a nice interval $C$. The transfer map to $C$ induced by $f$, 
$$T: \bigcup_{n \ge 0} f^{-n} (C) \mapsto C,$$
is defined as $T(x)=f^{\tau(x)}(x)$, where  $\tau:  \bigcup_{n \ge 0} f^{-n} (C) \mapsto \field{N}$, is the transfer time to $C$, that is the smallest non-negative integer $n$ such that $f^n(x) \in C$. 
\end{definition}

\begin{prop} \label{prop37}(Proposition 3.7 in \cite{Win2})
Assume that $f$ has no periodic attractors and that $S_f <0$. Let $T$ be the transfer map of $f$ to a   nice interval $C$. Then the complement if the domain of $T$ is a compact, $f$-invariant and hyperbolic set.
\end{prop}
\begin{proof}
Let $U={\rm dom} \ T$ and $\Gamma=[0,1] \setminus U$. Since $U$ is open, $\Gamma$ is closed, and, being bounded, is compact. By definition $f^{-1}(U) \subset U \implies f(\Gamma) \subset \Gamma$. 

$\Gamma$ is the set of points $x$ such that $f^n(x) \notin C$ for all $n \ge 0$. Since $S_f<0$, $f$ does not have non-hyperbolic periodic points (cf \cite{Mis}, Theorem 1.3), and, by assumption, $f$ has no periodic attracting orbits, so $\Gamma$  is hyperbolic (cf \cite{dMvS}, Lemma III.2.1).
\end{proof}

\medskip

Since a compact, invariant, hyperbolic set for a $C^{1+\alpha}$ function has zero Lebesgue measure  (cf. \cite{dMvS}, Theorem III.2.6), we have

\begin{cor}
$[0,1] \setminus {\rm dom} \  T$  has zero Lebesgue measure.
\end{cor}

The last result that we will require for the proof of unique ergodicity is the following Theorem due  Gambaudo and Martens (cf \cite{GaMa}).
\begin{thm}\label{GMtheorem}
If $f$ is infinitely renormalizable (of any combinatorial type) with a Cantor attractor $\Lambda$, then $\Lambda$ supports one or two ergodic invariant probability measures.

If the combinatorial type of $f$ is bounded, then $\Lambda$ is uniquely ergodic. 
\end{thm}

We can now present a proof of Theorem B, which is identical to Theorem $5.3$ from \cite{Win2}.

{\it Proof of Theorem B.} Let $L_n$ and $R_n$ denote the left and right half intervals for the $n$-th first return map, and let $i_n$ and $j_n$ be the return times for the corresponding intervals. Set $\Lambda_0=[0,1]$, and
$$\Lambda_n=\bigcup_{i=0}^{i_n-1}\overline{f^i(L_n)} \cup \bigcup_{i=0}^{j_n-1}\overline{f^j(R_n)},$$
where $\overline { \cdot}$ stands for the closure of a set. Components of $\Lambda_n$ are called intervals of generation $n$ and components of $\Lambda_{n-1} \setminus \Lambda_n$ are called gaps of generation $n$. Let $J \subset I$ be intervals of generations $n+1$ and $n$, respectively, and let $G \subset I$ be a gap of generation $n+1$.  Take the $\cL^0$ closure of the set $\{\cR^n[f]\}$. Since $\{\cR^n[f] \}$ is compact in $\cL^0$, the infimum and supremum of $|J|/|I|$ and $|G|/|I|$ over $I$, $J$ and $G$ of the corresponding generation are bounded away from $0$ and $1$. Otherwise, there would be an infinitely renormalizable map in $\cL^0$ with $|J|=0$ or $|J|=|I|$ ( $|G|=0$ or $|G|=|I|$). This is impossible, since this would imply that, for that map, one (or both) of $L_n$ or $R_n$ is of zero length, which contradicts renormalizability. Therefore, there exist constants $\mu>0$ and $\lambda<1$, such that
\begin{equation}\label{Cbounds}
\mu < {|J| \over |I|} < \lambda, \quad \mu < {|G| \over |I|} < \lambda.
\end{equation}

Next, $\Lambda \subset \cap \Lambda_n$, since the critical values are contained in the closure of $f(L_n) \cap f(R_n)$ for each $n$. From  $(\ref{Cbounds})$, $|\Lambda_{n+1}| \le \lambda |\Lambda_n|$, therefore, the lengths of intervals of generation $n$ tend to zero, and  $\Lambda=\cap \Lambda_n$.

A standard argument demonstrates that $\Lambda$ is a Cantor set of measure zero (since $\lambda<1$), of Hausdorff dimension in $(0,1)$. 

Next, we prove that almost all points are attracted to $\Lambda$. Let $T_n$ denote the transfer map to the $n$-th interval $C_n=L_n \cup R_n \cup \{0\}$. By Proposition $\ref{prop37}$ the domain of $T_n$ has full Lebesgue measure for every $n$, and, therefore, a.e. point enters $C_n$ for every $n$.

Finally, the unique ergodicity follows from Theorem $\ref{GMtheorem}$.
\qed

\bigskip

\section{Acknowledgments}
I would like to extend my gratitude to the families of Hernan Franco and Mariana Amadei, and  Pedro Beltran and Mariana Franco, for their hospitality during the preparation of this paper.



\begin{thebibliography}{99}
\bibitem{CCT} P. Collet, P. Coullet and C. Tresser, {\it Scenarios under  constraint}, J. Physique Lett. {\bf 46}(4) (1985), 143--147.
\bibitem{GaMa} J-M. Gambaudo and M. Martens, {\it Algebraic topology for minimal Cantor sets}, Ann. Henri Poincar\'e. {\bf 7}(3) (2006), 423--446.
\bibitem{GG} D. Gaidashev and I. Gorbovicki, {\it Complex bounds for Lorenz maps},  preprint (2016).
\bibitem{GD} A. Granas and J. Dugundji, {\it Fixed point theory},  Springer Monographs  in Mathematics (2003), Springer-Verlag, New-York.
\bibitem{GaiWin} D. Gaidashev and B. Winckler, {\it Existence of  a Lorenz renormalization fixed point of an arbitrary critical order},  preprint.
\bibitem{GW} J. Guckenheimer and R. F. Williams, {\it Structural stability of the Lorenz attractors},  Publ. Math. IHES {\bf 50} (1979), 59--72.
\bibitem{HS} J. H. Hubbard, C. T. Sparrow, {\it The classification of topologically expansive Lorenz maps},  Comm. Pure Appl. Math.  {\bf XLIII} (1990), 431--443.
\bibitem{dMvS} W. de Melo and S. van Strien, {\it One-dimensional dynamics}, {\bf 25} Ergbnisse der Matehematikund ihrer Grenzgebiete (3)  Ergod. Theor. and Dyn. Sys. {\bf 21}(3) (1993), Springer-Verlag, Berlin.
\bibitem{MM} M. Martens and W. de Melo, {\it Universal models for Lorenz maps},  Ergod. Theor. and Dyn. Sys. {\bf 21}(3) (2001), 883--860.
\bibitem{Win2} M. Martens and B. Winckler, {\it On the Hyperbolicity of Lorenz Renormalization}, Comm. Math. Phys., {\bf 325:1} (2014), 185-257.
\bibitem{Mis} M. Misiurewicz, {\it Absolutely continuous measures for certain maps of an interval}, Inst. Hautes \'Etudes Sci. Publ. Math. {\bf 53}, (1981), 17--51.
\bibitem{Lor} E. N. Lorenz, {\it Deterministic non-periodic flow},  J. Atmos. Sci {\bf 20} (1963), 130--141.
\bibitem{Wil} R. F. Williams, {\it The structure of the Lorenz attractors},  Publ. Math. IHES {\bf 50} (1979), 73--79.
\bibitem{Win1} B. Winckler, {\it A renormalization fixed point for Lorenz map},  Nonlinearity {\bf 23}(6) (2010), 1291--1303.


\end{thebibliography}
\end{document}